\documentclass[preprint,authoryear,1p,11pt]{elsarticle}

\usepackage{graphicx}
\DeclareGraphicsExtensions{.pdf,.png,.jpg}

\usepackage{amsmath}
\usepackage{amssymb}

\usepackage{multirow,dcolumn}

\usepackage{vmargin}
\setpapersize{USletter}
\setmarginsrb{2.5cm}{2.5cm}{2.5cm}{2.5cm}{0cm}{0cm}{0.5cm}{1cm}

\usepackage[mathcal,mathscr]{euscript}
\usepackage{xparse}
\usepackage{enumitem}
\usepackage{calc}

\usepackage{setspace}

\usepackage{tikz}
\usetikzlibrary{decorations.text,patterns,fadings,shapes.arrows,snakes}

\usetikzlibrary{external}
\tikzset{
    png export/.style={
        external/system call/.add={}{; convert -density 300 -transparent white "\image.pdf" "\image.png"},
        /pgf/images/external info,
        /pgf/images/include external/.code={
            \includegraphics[width=\pgfexternalwidth,height=\pgfexternalheight]{##1.png}
        },
    }
}

\usepackage{lineno}

\journal{Tree Physiology}

\newcommand{\leavethisout}[1]{}
\newcommand{\myee}[2]{\mbox{${#1}\times 10^{#2}$}}
\newcommand{\Afvwall}{\mathcal{A}_{\!f\!v}}
\newcommand{\Aroot}{\mathcal{A}_r}
\newcommand{\Atree}{\mathcal{A}_{\text{\emph{tree}}}}
\newcommand{\degC}{\mbox{${}^\circ$C}}

\newcommand{\myfiglett}[1]{{\bfseries\sffamily #1}}
\newcommand{\Creflect}{\mathcal{C}_r}
\newcommand{\Cin}{\mathcal{C}_{\text{\emph{r,in}}}}
\newcommand{\Cout}{\mathcal{C}_{\text{\emph{r,out}}}}
\newcommand{\Henry}{\mathcal{H}}
\newcommand{\heartfrac}{\theta}
\newcommand{\pSoil}{p_{\text{\emph{soil}}}}

\newcommand{\pIce}{p_{\text{\emph{ice}}}}
\newcommand{\Uroot}{U_r}
\newcommand{\Rgas}{\mathscr{R}}
\newcommand{\Rtree}{R_{\text{\rmfamily\emph{tree}}}}
\newcommand{\Rsap}{R_{\text{\rmfamily\emph{sap}}}}
\newcommand{\Tcrit}{T_c}
\newcommand{\TcSap}{T_{\text{\emph{c,sap}}}}
\newcommand{\Tamb}{T_a}

\newcommand{\tyreecond}{\mathscr{L}}
\newcommand{\Tmacro}{T}
\newcommand{\Emacro}{E}
\newcommand{\Xmacro}{\mathcal{X}}
\newcommand{\Tmicro}{\mathcal{T}}
\newcommand{\Emicro}{\mathcal{E}}
\newcommand{\Ymicro}{\mathcal{Y}}
\newcommand{\Yfast}{\Ymicro^1}
\newcommand{\Yslow}{\Ymicro^2}
\newcommand{\nderiv}[1]{\partial_n #1}
\newcommand{\yderiv}[1]{\partial_y #1}
\newcommand{\myrho}{\varrho}

\usepackage{caption}
\usepackage{ifthen}

\newboolean{@IsTikzPlotsOnly}
\setboolean{@IsTikzPlotsOnly}{false}
\setboolean{@IsTikzPlotsOnly}{true}

\newboolean{@IsProofs} 
\setboolean{@IsProofs}{false}
\setboolean{@IsProofs}{true}

\ifthenelse{\boolean{@IsTikzPlotsOnly}}{%
}{}

\usepackage[normalem]{ulem}

\definecolor{myBrickRed}{rgb}{0.80,0.20,0.20}
\definecolor{myTurquoise}{rgb}{0.25,0.68,0.61}
\definecolor{myCerulean}{rgb}{0.16,0.32,0.75}
\definecolor{myGray}{rgb}{0.7,0.7,0.7}
\newcommand{\mymod}[1]{{\color{myBrickRed} #1}}
\newcommand{\mydel}[1]{{\color{myGray}\sout{#1}}}
\renewcommand{\mymod}[1]{#1}
\renewcommand{\mydel}[1]{}

\begin{document}

\begin{frontmatter}



  \title{Experimental and computational comparison of freeze--thaw
    induced pressure generation in red and sugar maple}
  


\author[ubc]{Maryam Zarrinderakht}
\address[ubc]{Department of Earth, Ocean and Atmospheric Sciences,
  University of British Columbia, 2207 Main Mall, Vancouver, BC,
  V6T~1Z4, Canada} 

\author[comsysto]{Isabell Konrad}
\address[comsysto]{Comsysto Reply GmbH, Tumblingerstra{\ss}e 23, 80337
  Munich, Germany} 

\author[pmrc]{Timothy R. Wilmot}

\author[pmrc]{Timothy D. Perkins}

\author[pmrc]{\mbox{}\\ Abby K. van~den~Berg}
\address[pmrc]{Proctor Maple Research Center, University of
    Vermont, 58 Harvey Road, Underhill, Vermont, 05489, USA} 


\author[sfu]{John M.\ Stockie\corref{jms}}
\address[sfu]{Department of Mathematics, Simon Fraser University,
  8888 University Drive, Burnaby, BC, V5A~1S6, Canada}
\cortext[jms]{Corresponding author (jstockie@sfu.ca)} 

\vspace*{-2cm}

\begin{abstract}
%
  Sap exudation is the process whereby trees such as sugar
  (\textit{Acer saccharum}) and red maple (\textit{Acer rubrum})
  generate unusually high positive stem pressure in response to repeated
  cycles of freeze and thaw. This elevated xylem pressure permits the
  sap to be harvested over a period of several weeks and hence is
  a major factor in the viability of the maple syrup industry. The
  extensive literature on sap exudation documents competing
  hypotheses regarding the physical and biological mechanisms that drive
  positive pressure generation in maple, but to date relatively little
  effort has been expended on devising mathematical models for
  the exudation process. In this paper, we utilize an existing model of
  Graf et al.\ [J. Roy. Soc. Interface 12:20150665, 2015] that describes
  heat and mass transport within the multiphase gas--liquid--ice mixture
  in the porous xylem tissue.  The model captures the inherent
  multiscale nature of xylem transport by including phase change and
  osmotic transport in wood cells on the microscale, which is
  coupled to heat transport through the tree stem on the macroscale.  A
  parametric study based on simulations with synthetic temperature data
  identifies the model parameters that have greatest impact
  on stem pressure build-up.  Measured daily temperature fluctuations
  are then used as model inputs and the resulting simulated pressures
  are compared directly with experimental measurements taken from mature
  red and sugar maple stems during the sap harvest season. The results
  demonstrate that our multiscale freeze--thaw model reproduces
  realistic exudation behavior, thereby providing novel insights into
  the specific physical mechanisms that dominate positive pressure
  generation in maple trees. 
  %
\end{abstract}


\end{frontmatter}



\section{Introduction}
\label{sec:intro}


Sugar maple \emph{(Acer saccharum)}, red maple \emph{(Acer rubrum)}, and
several other \emph{Acer} species have a remarkable ability to generate
positive xylem (or sapwood) pressure in response to freeze--thaw cycles,
during a season when the tree is leafless and mostly dormant.  The
resulting exudation pressure can persist on and off for weeks or even
months, which allows maple sap to be harvested in sufficient quantities
that it is an economically viable agricultural product in northeastern
North America. A few other tree species such as black walnut
\emph{(Juglans nigra)}, butternut \emph{(Juglans cinerea)} and white
birch \emph{(Betula papyrifera)} are likewise capable of generating
elevated positive stem pressures, but none to such a high degree as
maple. It is well-known that the presence of dissolved sugar (2\%\ or
more by mass) in the sap of these species plays an important role in the
accumulation of stem pressure during the freeze--thaw process
\citep{marvin-1968}, but the precise causes of sap exudation have proven
difficult to pinpoint. Indeed, exudation has been studied extensively
for over 150 years, during which time researchers have attributed it to
a wide variety of physical and biological mechanisms, including osmosis
\citep{wiegand-1906}; thermal expansion of gas, water and wood
\citep{sachs-1860, merwin-lyon-1909, marvin-1949}; \mymod{active
  processes in living cells \citep{johnson-1945, marvin-1958}; and
  enhanced sap uptake due to ice crystal formation
  \citep{stevens-eggert-1945}.}

A major advance in the understanding of exudation was achieved by
\citet{milburn-omalley-1984}, whose experimental observations inspired
them to propose a purely physical hypothesis for maple sap exudation
based on a freeze--thaw mechanism that incorporates distinctive features
of the cellular structure of maple wood. More specifically, their model
focused on two classes of cells in the maple xylem: libriform fibers
that are filled with gas and are considered to play mainly a structural
role; and vessels and tracheids that are mostly liquid-filled and
constitute the primary pathways for sap transport.
\citeauthor{milburn-omalley-1984}'s breakthrough came from recognizing
that during a freezing cycle, \mymod{sap is drawn into the normally
  gas-filled fibers in the form of ice crystals that freeze on the
  interior surface of fiber walls thereby compressing the gas trapped
  within. This freeze-induced sap uptake is driven by an ice-water
  surface tension, which was already suggested by
  \citet{stevens-eggert-1945} to operate in sugar maple trees, and is
  analogous to the phenomenon of ``cryostatic suction'' that has been
  used to explain frost heave and ice lens formation in soils
  \citep{fowler-krantz-1994}.}  During a subsequent thaw, the stored
pressure is then released into the vessel sap.  This model explains the
variations in xylem pressure due to transfer of liquid during a sequence
of freeze--thaw events,
however it remains incomplete because it fails to explain the essential
role of dissolved sugars which are known to play a major role in
exudation.

This gap in understanding was addressed by several authors
\citep{tyree-1983, johnson-tyree-dixon-1987, johnson-tyree-1992} who
recognized that one consequence of the sugar content in sap is to
suppress the natural tendency in gas--liquid suspensions for gas bubbles
to dissolve at high pressure. \citet{tyree-1995} then proposed a
modification of the Milburn-O'Malley model that incorporates two
additional physical effects: expansion, contraction and dissolution of
gas bubbles in response to pressure variations in the suspending fluid;
and existence of an osmotic potential due to differences in solute
concentration engendered by the selectively permeable nature of certain
cell walls in maple xylem. In particular, \citeauthor{tyree-1995} argued
that the lignified cell wall separating the fiber and vessel
permits water and small solutes to pass but impedes larger molecules
like sucrose that make up the bulk of solutes in maple sap. This gives
rise to a significant osmotic potential between vessels (containing
high-sucrose sap) and surrounding fibers (containing pure water).
\citeauthor{tyree-1995} then demonstrated that this osmotic potential is
sufficient to stabilize gas bubbles over long enough time periods to
sustain a realistic exudation pressure, and the role of the selectively
permeable fiber/vessel walls in osmosis was later confirmed
experimentally by \citet{cirelli-etal-2008}. An extensive review of the
current literature on positive stem pressure and its causes can be found
in the paper by \citet{schenk-jansen-holtta-2021}.

Despite these advances in understanding of physical mechanisms behind
pressure generation in sugar maple and related species, significant
controversy remains in the literature over the causes of exudation.  For
example, \citet{ameglio-etal-2001} argued that ``no existing single
model explains all of the winter xylem pressure data,'' supported by
experiments suggesting that biological processes in living wood cells
are necessary for exudation. These differences in opinion have been
exacerbated by the lack of a mathematical description for the exudation
process. Indeed, around the time that \citeauthor{milburn-omalley-1984}
developed their model, \citet{tyree-1983} commented that ``there is
insufficient quantitative information to set up a system of physical
equations to describe the model''.  Some efforts have since been made to
formulate equations governing certain aspects of the relevant physics,
such as the diffusion model for embolism recovery developed by
\citet{yang-tyree-1992} that captures the gas transport and dissolution
processes similar to that which occurs in exudation.  

However, it was only in a series of papers
\citep{ceseri-stockie-2013, ceseri-stockie-2014, graf-stockie-2014,
  graf-ceseri-stockie-2015} that a concerted attempt was made to devise
a complete mathematical model of the essential physical mechanism behind
the freeze--thaw process of \citet{milburn-omalley-1984}, modified to
include the influence of gas bubbles and osmotic pressure
\citep{tyree-1995, cirelli-etal-2008}. The model is based on a system of
nonlinear diffusion equations for the cellular scale freeze--thaw
process that incorporate sap phase change, compression and dissolution
of gas bubbles, and osmosis \citep{ceseri-stockie-2013,
  ceseri-stockie-2014}, and captures the pressure exchange between
fibers and vessels during a single freeze--thaw cycle. This
cellular-level model was then coupled with a macroscale heat transport
equation obtained by a multi-scale averaging process known as periodic
homogenization \citep{graf-ceseri-stockie-2015}.
Numerical simulations yielded realistic pressure oscillations as well as
a build-up in exudation pressure over multiple freeze--thaw
cycles. These homogenized model simulations were then compared with
laboratory experiments on black walnut trees subjected to an imposed
periodic variation in temperature \citep{ameglio-etal-2001} and
exhibited excellent qualitative agreement.  A recent paper by
\citet{reid-driller-watson-2020} compared experimental measurements of
maple stem temperatures with solutions to a simpler but closely-related
model of heat transfer only, further demonstrating the effectiveness of this
class of models in predicting realistic stem temperatures.

Our main aim in this paper is to further validate the
model of \citet{graf-ceseri-stockie-2015} using experimental
measurements consisting of temperature and pressure time-series sampled
over a period of roughly 40 days
from red and sugar maple trees located at the University of Vermont
Proctor Maple Research Center in Underhill, Vermont, USA.  Having such a
highly-resolved dataset provides an ideal opportunity to validate the
time-dependent model for maple sap exudation under actual field
conditions.

After describing the experimental procedures used to collect field
measurements on maple trees in Section~\ref{sec:experiments} we present
the multiscale sap exudation model in Section~\ref{sec:model}, which
states the primary simplifying assumptions and supplies sufficient
details to elucidate the essential physical processes underlying the
model.
The comparison of experimental and numerical results in
Sections~\ref{sec:results}--\ref{sec:discussion} begins with a parameter
sensitivity analysis for a model tree subjected to a synthetic (smoothly
varying) ambient temperature, which aims to determine best estimates for
parameters used in our simulations of real maple trees.  We then perform
a sequence of simulations using actual data that feature rapid
variations in temperature and pressure, and significant levels
of measurement noise. A detailed comparison is then drawn
between stem pressures from simulations and experiments in order to
validate our hypothesis that a multiscale model based solely on physical
transport processes is capable of reproducing realistic exudation
pressures in maple trees subjected to repeated cycles of freeze and
thaw.

\section{Materials and Methods: Field Experiments}
\label{sec:experiments}

\begin{figure}[tbhp]
  \centering\small
  \ifthenelse{\boolean{@IsProofs}}%
  {\includegraphics[width=0.5\textwidth]{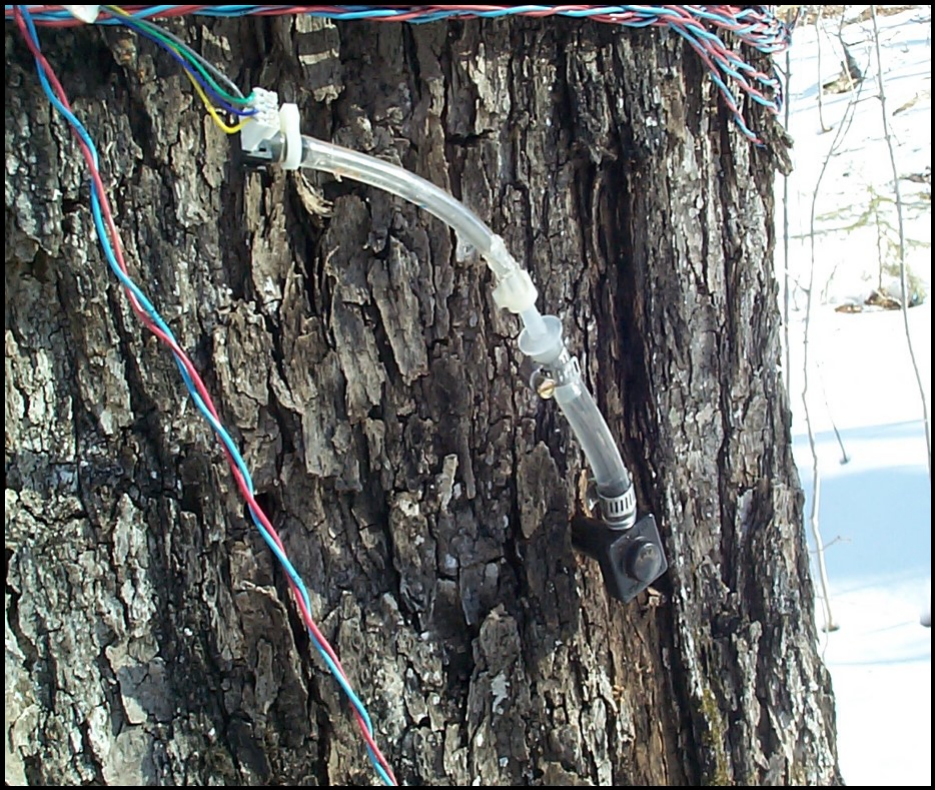}}%
  {\includegraphics[width=0.5\textwidth]{TPsensorPhoto}}
  \caption{Photograph of a taphole containing a black nylon spout that
    is connected by plastic tubing to an Omega PX-26-030GV pressure
    sensor.  Additional wires leading to thermocouples are also shown.}
  \label{fig:sensor-photo}
\end{figure}

Measurements were taken during the spring sap flow season over the years
2005--2010 in several mature sugar maple and red maple trees ranging in
size from 14 to 61~cm DBH (diameter at breast height). The trees were
located 180~m within the experimental sugarbush to the southeast of the
laboratory of the University of Vermont Proctor Maple Research Center in
Underhill, Vermont, at approximately 425~m elevation.  Type-T (24-gauge)
thermocouples were used to measure the following: air temperature at
heights of 1.2~m and 16.5~m; branch temperature measured in the center
of a hole drilled in a canopy branch at 16.5~m height, and within the
tapping zone at 1.2~m height and 5~cm depth; and soil temperature at
depths of 0~cm (surface) and 30~cm.  Pressure measurements were made by
drilling a standard 1.1~cm diameter, 5~cm deep taphole and inserting a
black nylon maple spout in an inverted position.  A tube leading from
the spout was filled with water and connected to an Omega PX-26-030GV
pressure sensor (with $\pm$1\% accuracy).  Both thermocouples and
pressure sensors were wired to a Campbell Scientific 21X datalogger and
measurements were stored at 15-minute intervals.  A photograph showing
the taphole and wired connections is given in
Figure~\ref{fig:sensor-photo}.  Stem temperature and pressure were
measured on both north and south sides of the tree and were often very
divergent.  Sap flow was measured from a taphole in a nearby tree of
similar size using a standard maple spout connected via tubing to a
collection chamber.  \mydel{A pressure sensor in the bottom was used to
  convert pressure to depth, so that exuded sap volume and flow rate
  could be determined for each interval.}  Data were downloaded from the
datalogger as comma-separated values and stored in Excel files for
analysis.

\leavethisout{
\begin{itemize}
\item Six red maples: diameter at breast height (DBH) $=$
  42.7, 28.4, 23.9, 18.5, 14.0, 25.7 cm
\item Two sugar maples: both DBH $\approx$ 61 cm
\end{itemize}
}

\section{Materials and Methods: Mathematical Model}
\label{sec:model}

\subsection{Background on xylem structure in maple}
\label{sec:xylem}

\begin{figure}[tbhp]
  \centering\small
  \ifthenelse{\boolean{@IsProofs}}%
  {\includegraphics[width=\textwidth]{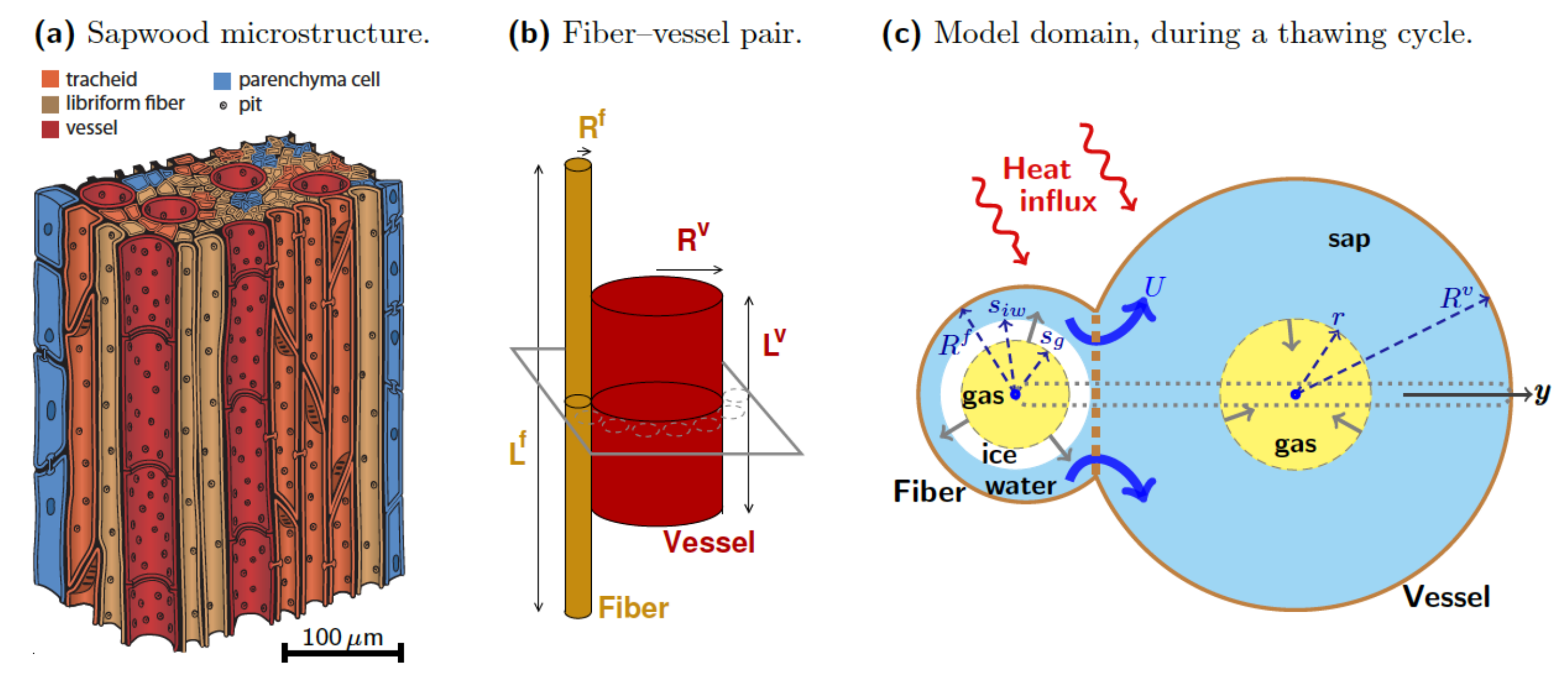}}%
  {\begin{tabular}{lclcl}
    \myfiglett{(a)} Sapwood microstructure. &~& 
    \myfiglett{(b)} Fiber--vessel pair. &~& 
    \myfiglett{(c)} Model domain, during a thawing cycle.\\[0.1cm]
    \ifthenelse{\boolean{@IsTikzPlotsOnly}}{%
    \begin{tikzpicture}
      \put(0,0) {\includegraphics[width=0.25\linewidth]{Images/trees/cells/XylemFinal}}
      \node at (4.1,6.45) {\mbox{}};
      \draw[ultra thin] (0,0) rectangle (0,0); 
      \draw[-, ultra thick] (2.8,0) -- (4.1,0);
      \draw[-, ultra thick] (2.8,-0.1) -- (2.8,0.1);
      \draw[-, ultra thick] (4.1,-0.1) -- (4.1,0.1);
      \node at (3.45,0.15) {\scriptsize\sffamily 100\,$\mu$m};
    \end{tikzpicture}
    }{\includegraphics[width=0.25\linewidth]{tikzfigure2a}}
    && 
    \raisebox{0.5cm}{\includegraphics[width=0.2\linewidth]{jrsint15-sapflow/fvvertical}}
    &&
    \ifthenelse{\boolean{@IsTikzPlotsOnly}}{%
    \begin{tikzpicture}[scale=1.2]
      \sffamily\bfseries
      \filldraw[color=brown, fill=cyan!40!white, very thick] (1,2) circle (1); 
      \filldraw[color=brown, fill=cyan!40!white, very thick] (3.5,2) circle (2); 
      \filldraw[color=cyan!40!white, fill=cyan!40!white] (1.5,1.3) rectangle ++(1,1.45);
      \filldraw[color=cyan!40!white, fill=white] (0.9,2) circle (0.7); 
      \filldraw[color=yellow, fill=yellow!80!white] (0.9,2) circle (0.5); 
      \draw[color=gray,very thin,densely dashed] (0.9,2) circle (0.5);
      \filldraw[color=yellow, fill=yellow!80!white] (3.5,2) circle (0.7); 
      \draw[color=gray,very thin,densely dashed] (3.5,2) circle (0.7);
      \node at (4.0,3.4) {{\footnotesize sap}}; 
      \node at (3.5,1.5) {{\footnotesize gas}};
      \node at (0.6,1.95) {{\footnotesize gas}};
      \node at (0.75,1.45) {{\footnotesize ice}};
      \node at (0.95,1.15) {{\footnotesize water}};
      \draw[color=blue!50!darkgray, ->, densely dashed, thick] (0.9,2) -- (1.2,2.4);
      \draw[color=blue!50!darkgray, ->, densely dashed, thick] (0.9,2) -- (0.8,2.7);
      \draw[color=blue!50!darkgray, ->, densely dashed, thick] (0.9,2) -- (0.4,2.8);
      \node[color=blue!60!black] at (1.25, 2.50) {\footnotesize $s_g$};
      \node[color=blue!60!black] at (0.83, 2.80) {\footnotesize $s_{iw}$};
      \node[color=blue!60!black] at (0.35, 2.49) {\footnotesize $R^f$};
      \draw[color=blue!50!darkgray, ->, densely dashed, thick] (3.5,2) -- (3.9,2.6);
      \draw[color=blue!50!darkgray, ->, densely dashed, thick] (3.5,2) -- (5.3,2.9);
      \node[color=blue!60!black] at (3.9,2.7) {\footnotesize $r$};
      \node[color=blue!60!black] at (5.0,2.9) {\footnotesize $R^v$};
      \draw[color=gray, ->, very thick] (1.22,1.62) -- (1.41,1.38);
      \draw[color=gray, ->, very thick] (0.46,1.76) -- (0.20,1.60);
      \draw[color=gray, ->, very thick] (1.02,2.48) -- (1.11,2.77);
      \draw[color=gray, ->, very thick] (4.11,1.65) -- (3.84,1.80);
      \draw[color=gray, ->, very thick] (2.84,1.76) -- (3.12,1.86);
      \draw[color=gray, ->, very thick] (3.44,2.70) -- (3.47,2.40);
      \coordinate (S) at (1.4,2.7);
      \coordinate (E) at (2.1,2.9);
      \draw[color=blue, ->, line width=0.1cm, opacity=0.8] (S) to [out=-70,in=-110,distance=0.4cm] (E);
      \coordinate (SS) at (1.4,1.2);
      \coordinate (EE) at (2.1,1.0);
      \draw[color=blue, ->, line width=0.1cm, opacity=0.8] (SS) to [out=70,in=110,distance=0.4cm] (EE);
      \node[color=blue] at (2.2,3.0) {\footnotesize $U$};
      \draw[color=gray, dotted, very thick, rounded corners] (0.9,1.9) rectangle ++(4.6,0.2);    
      \draw[->, color=darkgray, thick] (4.5,2) -- (5.7,2);
      \node at (5.8,2) {$\pmb{y}$};
      \draw[-, color=blue, very thick] (0.9,2) circle (0.03);
      \draw[-, color=blue, very thick] (3.5,2) circle (0.03);
      \draw[-, dashed, line width=0.8mm, color=brown] (1.65,1.26) -- (1.65,2.73); 
      \node at (0.1,1.11) {Fiber};
      \node at (4.9,0.13) {Vessel};
      \draw[->, color=red!80!darkgray, snake=coil, segment aspect=0, very thick] (0.5,4.0) -- (1.0,3.2); 
      \draw[->, color=red!80!darkgray, snake=coil, segment aspect=0, very thick] (1.5,4.5) -- (2.0,3.7); 
      \node at (1.1,4.1) {\color{red!80!darkgray}Heat};
      \node at (1.3,3.8) {\color{red!80!darkgray}influx};
      \draw[-, color=white] (0.5,-0.5) -- (0.5,-0.5);
    \end{tikzpicture}
  }{\raisebox{0.25cm}{\includegraphics[width=0.45\linewidth]{tikzfigure2c}}}
  \end{tabular}}
  \caption{Sapwood microstructure and the idealized 2D model geometry.
    (a)~A microscopic cut-away view of the sapwood within a typical
    hardwood tree, depicting the vessels and (libriform) fibers that are
    central to the model.  Tracheids are connected hydraulically to
    neighboring vessels via paired pits, which is why they are ``lumped
    together'' in our model with vessels. Note that fiber walls also
    contain pits, but they are unpaired and hence unconnected to
    adjacent vessels or tracheids.  (b)~A single fiber--vessel pair
    showing the main geometric parameters. The horizontal cutting plane
    highlights the planar cross-section corresponding to the 2D model
    geometry in Figure~{\protect\ref{fig:cell-geometry}}c. The dashed
    circles represent the $N$ copies of the fiber that are incorporated
    into the equations through a simple multiplier $N$. (c)~The 2D model
    geometry depicting a thawing scenario. A thawing fiber of radius
    $R^f$ (containing nested layers of gas, ice and liquid water) is
    located adjacent to a thawed vessel of radius $R^v$ (containing gas
    and liquid sap). As the fiber ice layer thaws, the fiber gas bubble
    expands and forces melt-water through the porous wall into the
    vessel at a rate $U$, thereby compressing the vessel gas and
    increasing the vessel sap pressure.}
  \label{fig:cell-geometry}
\end{figure}

Before presenting the mathematical model for the freeze--thaw process
governing sap exudation, we begin by briefly summarizing the
physical and structural characteristics of maple xylem that play an
essential role in sap transport and exudation (more detail is provided
in \citet{tyree-zimmermann-2002}).  The xylem in hardwood trees such
as maple is built primarily of rigid and nearly cylindrical structures
that consist of the hollowed-out walls of dead wood cells. Sapwood has a
regular and quasi-periodic microstructure shown in
Figure~\ref{fig:cell-geometry}a that consists of regularly spaced
vessels interspersed with much more abundant tracheids and (libriform)
fibers having a smaller diameter. Vessels are the primary
water-conducting conduits in sapwood, and each is divided lengthwise
into vessel elements that are connected end-to-end via perforation
plates, thereby forming long capillary tubes. Vessel walls are
interspersed with cavities called pits that connect cells
hydraulically through pit membranes as long as the pits in adjacent cells
align or ``pair up''. Tracheids are intermediate in size between vessels
and fibers and are also connected through paired pits with neighboring
vessels and other tracheids. Since both vessels and tracheids are mostly
sap-filled, they serve as the primary conduits for xylem sap transport.
Fibers on the other hand are known to contain mostly
gas~\citep{milburn-omalley-1984} and their walls demonstrate a relative
lack of pitting. As a result, they are usually regarded to have
negligible impact on sap transport and instead serve a structural
function. The lignified fiber walls have nonetheless been found to
contain micropores (much smaller than those in pit membranes) which are
selectively permeable, in the sense that they allow water and small
solutes to pass but inhibit the passage of larger molecules like sucrose
\citep{cirelli-etal-2008}.

\subsection{Microscale model for the thawing process}
\label{sec:model-thawing}

We next describe our model for the essential processes governing sap
exudation on the microscopic (cellular) scale which include: heat
transport; phase change due to freezing and thawing; expansion,
contraction and dissolution of gas; and porous flow through paired pit
membranes and selectively-permeable micropores in the fiber--vessel
wall. Our ultimate goal is to model repeated freeze--thaw cycles in
response to ambient temperature variations as depicted in
Figure~\ref{fig:thawing-freezing}, where a thawing (or freezing) front
propagates into the tree stem and separates regions of the sapwood that
are frozen from those that are thawed.  The front actually comprises a thin
annular section of sapwood that exists in a mixed state containing
water/sap in both liquid and ice phases.  We focus our attention on
equations for the thawing process only and refer the interested reader
to \citet{graf-ceseri-stockie-2015} for the modifications required to
capture all other possible freeze--thaw stages.

\tikzfading[name=arrowfading, top color=transparent!0, bottom color=transparent!95]
\tikzset{arrowfill/.style={#1}}
\tikzset{arrowstyle/.style n args={3}{draw=#2,arrowfill={#3}, 
    single arrow, minimum height=#1, single arrow, 
    single arrow head extend=0.15cm}}
\NewDocumentCommand{\tikzfancyarrow}{O{2cm} O{darkgray} 
  O{top color=darkgray!20, bottom color=darkgray} m}{
  \tikz[baseline=-0.5ex]\node [arrowstyle={#1}{#2}{#3}] {#4};}

\begin{figure}[tbhp]
  \centering\sffamily\bfseries
  \ifthenelse{\boolean{@IsProofs}}%
  {\includegraphics[width=0.65\textwidth]{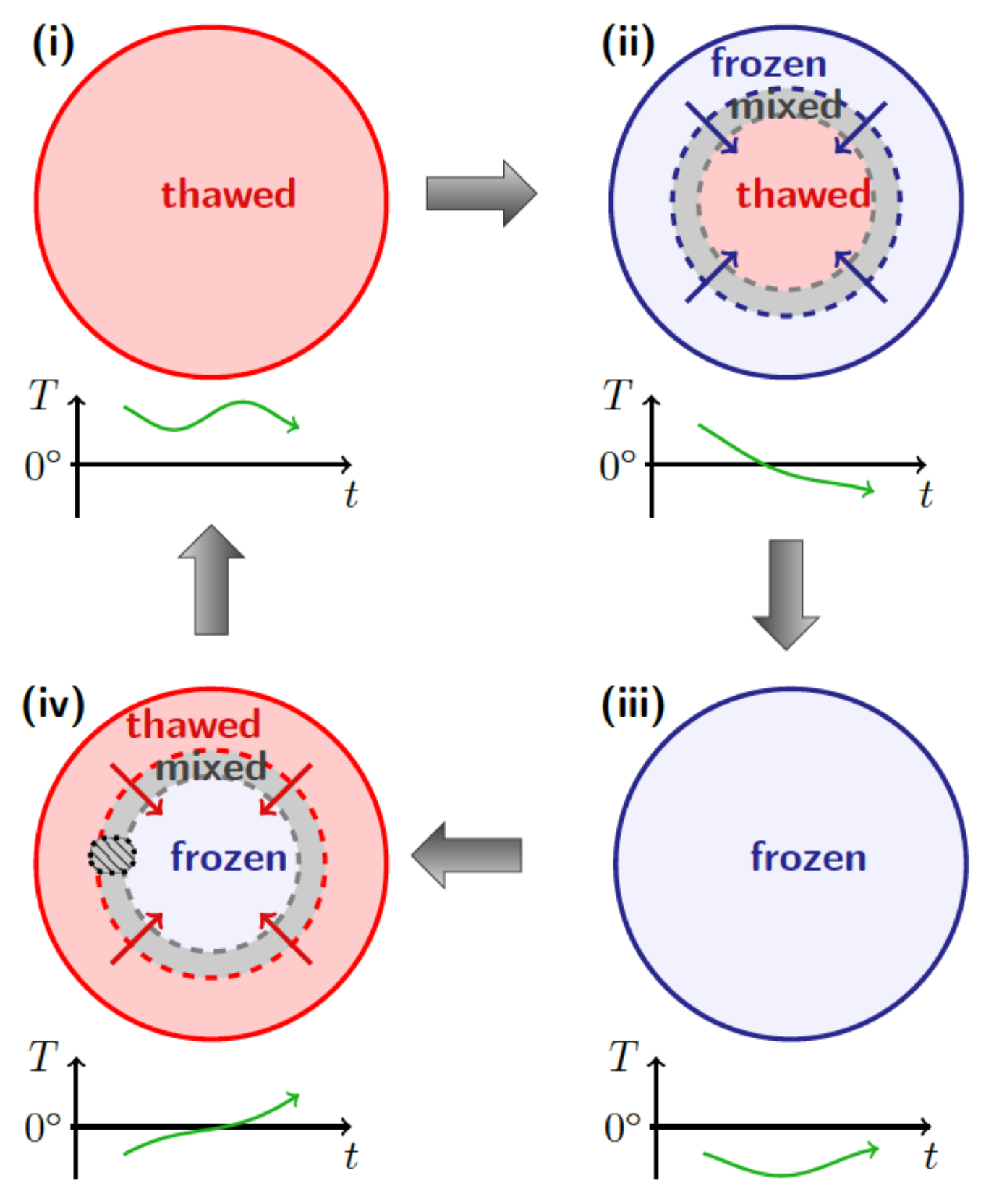}}%
  {%
  %
  \ifthenelse{\boolean{@IsTikzPlotsOnly}}{%
  \begin{tikzpicture}[scale=1.6]
    \node[color=black] at (0.1,1.9) {(i)};
    \filldraw[color=red, fill=red!20, very thick](1,1) circle (1);
    \node[align=center,color=red!80!darkgray]  at (1.1,1.05) {thawed};
    \draw [thick,->] (0.2,-0.5)--(1.8,-0.5) node [right, below]{$t$};
    \draw [thick,->] (0.23,-0.8)--(0.23,-0.1) node [left]{$T$};
    \node at (0.25,-0.5) [left] {$0^\circ$};
    \draw [->, thick, variable=\t, domain=0.5:1.5, samples=101, smooth, color=green!60!darkgray] 
       plot (\t, {-0.22-0.08*cos(8*\t r))});
    %
    \node at (2.5,1.05) {\tikzfancyarrow[1cm]{}};
    \node[rotate=90] at (1.0,-1.2) {\tikzfancyarrow[1cm]{}};
  \end{tikzpicture}
  }{\includegraphics[width=0.33\linewidth]{tikzfigure3i}}
  %
  \ifthenelse{\boolean{@IsTikzPlotsOnly}}{%
  \begin{tikzpicture}[scale=1.6]
    \node[color=black] at (0.1,1.9) {(ii)};
    \filldraw[color=blue!40!darkgray, fill=blue!5, very thick] (1,1) circle (1);
    \filldraw[color=blue!40!darkgray, fill=gray!40!white, very thick, dashed] (1,1) circle (0.65);
    \filldraw[color=gray, fill=red!20, very thick, dashed] (1,1) circle (0.50);
    \draw[color=blue!40!darkgray, ->, very thick] (1.57,1.57) -- (1.28,1.28);
    \draw[color=blue!40!darkgray, ->, very thick] (0.43,1.57) -- (0.72,1.28);
    \draw[color=blue!40!darkgray, ->, very thick] (0.43,0.43) -- (0.72,0.72);
    \draw[color=blue!40!darkgray, ->, very thick] (1.57,0.43) -- (1.28,0.72);
    \node[align=center,color=red!80!darkgray]  at (1.1,1.05) {thawed};
    \node[align=center,color=blue!40!darkgray] at (0.9,1.8) {frozen};
    \node[align=center,color=darkgray] at (1,1.56) {mixed};
    \draw [thick,->] (0.2,-0.5)--(1.8,-0.5) node [right, below]{$t$};
    \draw [thick,->] (0.23,-0.8)--(0.23,-0.1) node [left]{$T$};
    \node at (0.25,-0.5) [left] {$0^\circ$};
    \draw [->, thick, variable=\t, domain=0.5:1.5, samples=101, smooth, color=green!60!darkgray] 
       plot (\t, {-0.1-0.4*\t+0.05*sin(5*\t r)});
    \node[rotate=270] at (1.0,-1.2) {\tikzfancyarrow[1cm]{}};
  \end{tikzpicture}
  }{\includegraphics[width=0.23\linewidth]{tikzfigure3ii}}

  %
  \ifthenelse{\boolean{@IsTikzPlotsOnly}}{%
  \begin{tikzpicture}[scale=1.6]
    \node[color=black] at (0.1,1.9) {(iv)};
    \filldraw[color=red, fill=red!20, very thick](1,1) circle (1);
    \filldraw[color=red, fill=gray!40!white, very thick, dashed] (1,1) circle (0.65);
    \filldraw[color=gray, fill=blue!5, very thick, dashed] (1,1) circle (0.50);
    \draw[color=red!80!darkgray, ->, very thick] (1.57,1.57) -- (1.28,1.28);
    \draw[color=red!80!darkgray, ->, very thick] (0.43,1.57) -- (0.72,1.28);
    \draw[color=red!80!darkgray, ->, very thick] (0.43,0.43) -- (0.72,0.72);
    \draw[color=red!80!darkgray, ->, very thick] (1.57,0.43) -- (1.28,0.72);
    \node[align=center,color=blue!40!darkgray] at (1.1,1.05) {frozen};
    \node[align=center,color=red!80!darkgray] at (0.9,1.8) {thawed};
    \node[align=center,color=darkgray] at (1,1.56) {mixed};
    \draw[color=gray, rounded corners, fill=gray!40!white] 
      (0.31,0.95) rectangle ++(0.25,0.20); 
    \draw[dotted, color=black, very thick, rounded corners, fill=gray, 
      pattern color=darkgray, pattern=north west lines]
      (0.31,0.95) rectangle ++(0.25,0.20);    
    \draw [thick,->] (0.2,-0.5)--(1.8,-0.5) node [right, below]{$t$};
    \draw [thick,->] (0.23,-0.8)--(0.23,-0.1) node [left]{$T$};
    \node at (0.25,-0.5) [left] {$0^\circ$};
    \draw [->, thick, variable=\t, domain=0.5:1.5, samples=101, smooth, color=green!60!darkgray] 
       plot (\t, {-0.9+0.4*\t-0.05*cos(5*\t r))});
    \node[rotate=180] at (2.5,1.05) {\tikzfancyarrow[1cm]{}};
  \end{tikzpicture}
  }{\includegraphics[width=0.33\linewidth]{tikzfigure3iv}}
  %
  \ifthenelse{\boolean{@IsTikzPlotsOnly}}{%
  \begin{tikzpicture}[scale=1.6]
    \node[color=black] at (0.1,1.9) {(iii)};
    \filldraw[color=blue!40!darkgray, fill=blue!5, very thick] (1,1) circle (1);
    \node[align=center,color=blue!40!darkgray]  at (1.1,1.05) {frozen};
    \draw [thick,->] (0.2,-0.5)--(1.8,-0.5) node [right, below]{$t$};
    \draw [thick,->] (0.23,-0.8)--(0.23,-0.1) node [left]{$T$};
    \node at (0.25,-0.5) [left] {$0^\circ$};
    \draw [->, thick, variable=\t, domain=0.5:1.5, samples=101, smooth, color=green!60!darkgray] 
       plot (\t, {-0.7+0.08*sin(5*\t r))});
  \end{tikzpicture}
  }{\includegraphics[width=0.23\linewidth]{tikzfigure3iii}}}
  \caption{The freeze--thaw process within a circular tree stem cycles
    between four main stages
    (i\,$\rightarrow$\,ii\,$\rightarrow$\,iii\,$\rightarrow$\,iv\,$\rightarrow$\,i\,$\rightarrow$\,\dots)
    as ambient temperature $T$ cycles below and above the freezing
    point: (i) completely thawed (with $T>0$); (ii) partially frozen
    ($T\searrow 0$), with a freezing front advancing radially inward to
    the center of the stem; (iii) completely frozen ($T<0$); (iv)
    partially thawed ($T\nearrow 0$), with a thawing front advancing
    radially inward. The freezing/thawing fronts in (ii,iv) are thin
    annular regions (shaded in grey, and in reality much thinner than
    depicted here) wherein the liquid is in a ``mixed'' state; that is,
    the water in the fibers is freezing/frozen and the sap in the
    vessels thawing/thawed. The thawing front circled on the left of
    (iv) is magnified in Figure~\ref{fig:cell-geometry}c to the cellular
    scale, which depicts an individual vessel and an adjacent fiber in a
    partially thawed state.}
  \label{fig:thawing-freezing}
\end{figure}

In order that our model remain tractable, we make a number of
simplifying assumptions:
\begin{enumerate}[label=A\arabic*.,ref=A\arabic*]
\item Vertical \label{assume:2d} variations due to gravity and
  height-dependence are neglected, so that we can focus on a 2D
  horizontal cross-section through a circular tree stem. \mymod{We also
    assume a constant stem diameter because there was no noticeable
    change in the stem cross-section during the period that measurements
    were taken.}
  
\item At \label{assume:cells} the cellular level, we consider only the
  contribution of vessels and fibers to exudation. Tracheids
  are not treated separately but rather lumped together with
  vessels. Other sapwood components (such as ray cells) are ignored
  but are acknowledged as the source of sucrose in the vessel
  lumens.
  
\item Sapwood \label{assume:periodicity} has a uniform, periodic
  microstructure consisting of cylindrical vessels and fibers, both
  having constant lengths ($L^v$, $L^f$) and radii ($R^v$, $R^f$). Then
  we can reasonably restrict our attention to a horizontal slice through
  a single vessel element and adjacent fiber as pictured in
  Figure~\ref{fig:cell-geometry}b.
  
\item Each \label{assume:Nfibers} vessel is in contact with $N$ fibers
  on average (see Figure~\ref{fig:cell-geometry}b) so that the influence
  of multiple fibers can be incorporated by simply multiplying by a
  factor of $N$ the contribution from the single fiber being modelled.
  
\item Because \label{assume:radial-symmetry} the fiber and vessel
    cross-sections are circular (from assumption
    \ref{assume:periodicity}) we assume radial symmetry and take all
    microscopic variables to depend on a radial coordinate $y$ that
    passes through the line joining the centers of the fiber and
    vessel.
  
  \item Gas \label{assume:vessel-gas} is also present in the vessels,
    which is consistent with observations showing that exuding maple sap
    contains suspended gas bubbles \citep{wiegand-1906,
      marvin-greene-1959, perkins-vandenberg-2009}. \mymod{We are also
      guided by the observation of \citet{milburn-omalley-1984} that
      both sap and xylem tissue are essentially incompressible at the
      pressures typically experienced during maple sap
      exudation. Consequently, gas bubbles suspended in the vessel sap
      provide a convenient mechanism to facilitate pressure exchange
      between fiber and vessel.}  As the sap pressure periodically rises
    and falls through the freeze-thaw cycles, a portion of the gas may
    dissolve within the sap if the local pressure is high enough.
  
\item Gas, \label{assume:layers} liquid and ice within the fiber and
  vessel exist as distinct layers arranged as concentric annuli shown in
  Figure~\ref{fig:cell-geometry}c. This assumption is made for
  mathematical convenience since it permits the geometry to be
  captured by a single radial coordinate.  This is also consistent
  with \citet{milburn-omalley-1984} who conjectured that ice accumulates
  on the inner fiber wall due to cryostatic suction from previous
  freezing cycles to encase a central gas bubble, and that during a
  thawing cycle any liquid melt-water collects in a layer between the
  ice and the wall.
  When it comes to the vessels, experimental observations suggest that
  bubbles form asymmetrically at nucleation sites along the inside of
  the vessel lumen~\citep{brodersen-etal-2010,
    zwieniecki-holbrook-2000}; however, the precise location of these
  bubbles is unimportant for the purposes of our model because
  \mymod{the influence of the gas appears only as a volume fraction} in
  the model equations.
  
  
\item The \label{assume:constant-T} temperature in the gas and ice
  layers is assumed to be constant and equal to that of the adjacent
  liquid. Furthermore, during any freeze or thaw cycle the temperatures
  of the gas/ice layers in the fiber remain constant at \mymod{the
    critical temperature (or freezing point) of water, denoted
    $\Tcrit$}. This is justified because the fiber radius is roughly 6
  times smaller than in the vessel, and the thermal diffusivities for
  gas and ice ($1.9\times 10^{-5}$ and $1.5\times 10^{-6}$
  m${}^{2}$s${}^{-1}$) are at least an order of magnitude larger than
  that for water ($1.4\times 10^{-7}$ m${}^{2}$s${}^{-1}$)
  \citep{tyree-zimmermann-2002}.

\item Any \label{assume:permeability} liquid entering the fiber from the
  vessel must have a negligible sucrose content owing to the
  selective permeability of the fiber--vessel wall \citep{cirelli-etal-2008}.
  
\item Liquid \label{assume:soil-liquid} water is present within the soil
  even under freezing conditions and tree roots passively
  transport water throughout the entire freeze--thaw cycle, both of
  which are supported by observations \citep{marvin-1958,
    robitaille-boutin-lachance-1995, sorkin-2014}.  However, we do
  recognize that the upper few centimeters of soil may remain frozen.
\end{enumerate}

Based on the above assumptions, we now present the governing equations
for the thawing cycle of the exudation process during which vessel are
completely thawed, whereas the fibers are in a partially frozen state
(depicted in Figure~\ref{fig:cell-geometry}c and corresponding to stage
(iv) in Figure~\ref{fig:thawing-freezing}).  \mymod{The model equations
  are essentially the same as those presented in
  \citet{graf-ceseri-stockie-2015}, and complete details of the
  derivation can be found there and in \citet{ceseri-stockie-2013}}.
The physical state of the various phases (liquid, ice, gas) within a
given vessel and fiber pair can be described by the following six
time-dependent functions:
\begin{itemize}
\item[] $s_g(t)$: fiber gas bubble radius, measured from the center of
  the fiber, 
\item[] $s_{iw}(t)$: radius of the fiber ice-water interface, 
\item[] $r(t)$: vessel gas bubble radius, 
\item[] $U(t)$: total volume of melt-water that flows through the
  porous fiber--vessel wall, measured positive from fiber to vessel,
\item[] $\Uroot(t)$: total volume of water influx from the
  roots,
\item[] $\Tmicro(y,t)$: temperature in the vessel sap, which also depends on
  the radial coordinate $y$ with origin at the vessel center so that
  $r(t)\leqslant y \leqslant R^v$ (recall that temperature is taken to
  be constant elsewhere in the fiber and vessel).
\end{itemize}
We will now state equations for the first five unknowns, leaving the
microscale heat equation for $\Tmicro$ to the next section.  First of
all, an algebraic equation for conservation of volume can be derived
that relates the thickness of various layers within the fiber and
balances with any melt-water exiting into the vessel ($U$).  After
exploiting the circular symmetry in the fiber, this volume constraint is
differentiated in time to obtain
\begin{linenomath}\begin{gather}
  \partial_t s_g = - \frac{(\myrho_w -\myrho_i) s_{iw} \partial_t
    s_{iw}}{\myrho_i s_g} + \frac{\myrho_w \partial_t U}{2 \pi
    L^f \myrho_i s_g}, 
  \label{eq:sgas}
\end{gather}\end{linenomath}
where $\myrho_w$ and $\myrho_i$ are the water and ice densities
respectively. Note that this equation is an ordinary differential equation (ODE)
for the variable $s_g(t)$.  A similar volume conservation equation
governs the vessel gas bubble radius $r(t)$
\begin{linenomath}\begin{gather}
  \partial_t r = - \frac{N \partial_t U + \partial_t \Uroot}{2\pi L^v
    r},  
  \label{eq:rgas}
\end{gather}\end{linenomath}
which involves an additional term arising from root water $U$ drawn into
the xylem vessel network.  Note that the $\partial_t U$ term is
multiplied here by $N$, which is the average number of fibers connected
to each vessel. The dynamics of the ice-water interface $s_{iw}$ is
governed by a phase change process that obeys the Stefan condition
\begin{linenomath}\begin{gather}
  \partial_t s_{iw} = -\frac{k_w/\myrho_w}{E_w-E_i} \, \nderiv{\Tmicro} +
  \frac{\partial_t U}{2 \pi L^f s_{iw}},  
  \label{eq:siw-stefan}
\end{gather}\end{linenomath}
where $k_w$ is the thermal conductivity of water, $(E_w-E_i)$ is the
latent heat, and $\nderiv{\Tmicro} \equiv \nabla_y \Tmicro\cdot\vec{n}$
denotes the normal heat flux.

The two remaining equations are obtained by applying
Darcy's law  for flow in porous media.  The liquid flux through the
porous fiber--vessel wall obeys 
\begin{linenomath}\begin{gather}
  \partial_t U = - \frac{\tyreecond \Afvwall}{N} \left( p_w^v - p_w^f
    - C_s \Rgas \Tmicro(R^v,t) + \pIce \right)
  \label{eq:U-darcy}
\end{gather}\end{linenomath}
where $\Afvwall = 2 \pi R^v L^v$ is the surface area of the
fiber--vessel wall and $\tyreecond$ is its hydraulic conductivity.  The
term in parentheses represents a balance at the fiber--vessel wall
between four pressures: vessel and fiber liquid pressures, $p_w^v(t)$
and $p_w^f(t)$; osmotic pressure 
deriving from the sap sugar concentration $C_s$; and cryostatic suction
$\pIce$ which is an ice-water surface tension that is zero under thawing
conditions, but nonzero when the fiber is completely frozen and the
adjacent vessel contains liquid sap.  A second application of Darcy's
law at the root gives the volume flux of root water as
\begin{linenomath}\begin{gather}
    \partial_t \Uroot = - \tyreecond_r \Aroot \mymod{(p_w^v - \Psi_s)},
  %
  %
  %
  \label{eq:Uroot-darcy}
\end{gather}\end{linenomath}
where $\tyreecond_r$ is the root hydraulic conductivity and $\Aroot$ is
the root area (per vessel).  The pressure balance here is expressed in
terms of a difference betwen \mymod{water potentials in the soil
  ($\Psi_s$) and xylem ($p_w^v$, our vessel sap pressure).}  \mymod{With
  this in mind, the conductivity $\tyreecond_r$ should actually be
  viewed as an effective parameter that account for the various
  conductivities along the water path through roots, xylem and vessel.}

%
\leavethisout{%
  According to \citet{henzler-etal-1999}, the roots function
  as a partial check valve in the sense that aquaporin membranes
  controlling root water transport are more permeable to inflow than to
  outflow. This effect is incorporated through a reflection coefficient
  $\Creflect \in [0,1]$ that takes the value $\Creflect=\Cin=1$ for
  inflow (when $p_w^v \leqslant \pSoil$) and $\Creflect=\Cout\in[0,1)$
  for outflow (when $p_w^v > \pSoil$).
  The lower limit $\Cout=0$ corresponds to the case of no root outflow
  (which was assumed by \citet{graf-ceseri-stockie-2015}) whereas the
  upper limit $\Cout=1$ represents the symmetric case where root inflow and
  outflow experience equal resistance.  We propose using an intermediate
  value of $\Cout=0.2$, which is consistent with experiments of
  \citet{henzler-etal-1999} on legumes, and of
  \citet{steudle-peterson-1998} on woody plants (although a
  values as large as $\Cout=0.7$ have been suggested by
  \citet{tyree-etal-1994}).%
}

The only quantities in the above equations that remain to be specified
are the liquid pressures $p_w^v$ and $p_w^f$ within the vessel and
fiber.  In the vessel, the sap pressure satisfies
\begin{linenomath}\begin{gather}
  p_g^v(t) = \frac{\myrho_g^v(t) \Rgas \Tcrit}{M_g}
  - \frac{2\sigma_{gw}}{r(t)},  
  \label{eq:pressure1}
\end{gather}\end{linenomath}
which represents a balance between gas pressure (from the ideal gas law) and
surface tension at the gas bubble interface
(from the Young-Laplace equation).  Here, $\myrho_g^v(t)$ is the vessel
gas density and $\sigma_{gw}$ is the gas-liquid surface tension.  The
gas density is related to the volumes $V_w^v(t)$ and $V_g^v(t)$ of the
vessel sap and gas regions by
\begin{linenomath}\begin{gather}
  \myrho_g^v(t) = \left(\frac{V_g^v(0) + \Henry V_w^v(0)}{V_g^v(t) +
      \Henry V_w^v(t)}\right) \myrho_g^v(0),  
  \label{eq:henry1}
\end{gather}\end{linenomath}
which accounts for gas dissolving in the sap via terms involving Henry's
constant $\Henry$.  Finally, the cell volume quantities are determined
by simple geometric constraints
\begin{linenomath}\begin{gather}
  V_g^v(t) = \pi L^v r(t)^2 \qquad \text{and} \qquad 
  V_w^v(t) = \pi L^v \left((R^v)^2 - r(t)^2\right).
  \label{eq:volume1}
\end{gather}\end{linenomath}
An analogous set of equations govern the fiber water pressure $p_w^f(t)$
but are not included here.


Recall that the ODEs~\eqref{eq:sgas}--\eqref{eq:Uroot-darcy} along with
the algebraic constraints~\eqref{eq:pressure1}--\eqref{eq:volume1}
describe the dynamics within the fiber and vessel during a thawing cycle
only, when a fiber in the midst of thawing lies adjacent
to a vessel that is completely thawed.  There are five
additional cases corresponding to the various stages of freeze--thaw
within the fiber and vessel, and each case leads to a
modification of the governing equations detailed in
\citet{graf-ceseri-stockie-2015}.

\subsection{Homogenized two-scale model for heat transport}
\label{sec:model-heat}

\begin{figure}[tbhp]
  \centering\bfseries\sffamily\small
  \ifthenelse{\boolean{@IsProofs}}%
  {\includegraphics[width=0.95\textwidth]{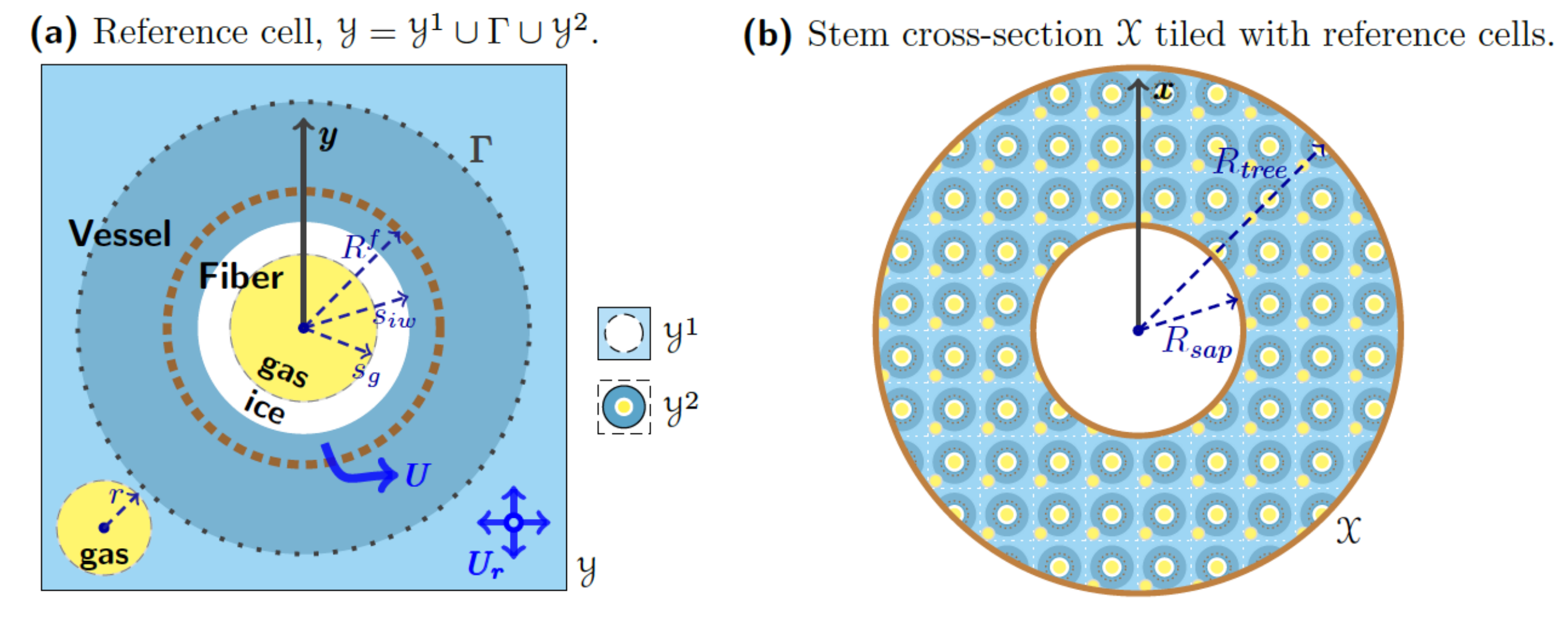}}%
  {%
  \ifthenelse{\boolean{@IsTikzPlotsOnly}}{%
  \begin{tikzpicture}
    \draw[ultra thin, fill=cyan!40!white] (0,0) rectangle (5,5); 
    \node at (2.6,5.3) {(a) {\rmfamily\normalfont Reference cell, $\Ymicro=\Yfast\cup\Gamma\cup\Yslow$.}};
    \draw[color=darkgray, loosely dotted, very thick, fill=cyan!60!gray!80!white] (2.5,2.5) circle(2.15);
    \node[color=darkgray] at (4.2,4.2) {{\small $\pmb{\Gamma}$}};
    \draw[color=brown!80!black, densely dashed, line width=0.8mm] (2.5,2.5) circle(1.3);
    \draw[color=white, fill=white] (2.5,2.5) circle(1.0);
    \node[rotate=-25] at (2.15,1.75) {{\footnotesize ice}};
    \draw[color=gray!30!white, fill=yellow!80!white] (2.5,2.5) circle(0.7);
    \draw[color=gray,very thin,densely dashed] (2.5,2.5) circle (0.7);
    \node[rotate=-25] at (2.3,2.05) {{\footnotesize gas}};
    \draw[color=blue!50!darkgray, ->, thick, densely dashed] (2.5,2.5) -- (3.42,3.42);
    \node[color=blue!60!black] at (3.05,3.30) {\footnotesize $R^f$};
    \draw[color=blue!50!darkgray, ->, thick, densely dashed] (2.5,2.5) -- (3.5,2.8);
    \node[color=blue!60!black] at (3.37,2.59) {\footnotesize $s_{iw}$};
    \draw[color=blue!50!darkgray, ->, thick, densely dashed] (2.5,2.5) -- (3.15,2.25);
    \node[color=blue!60!black] at (3.10,2.05) {\footnotesize $s_g$};
    \draw[color=gray!30!white, fill=yellow!80!white] (0.6,0.6) circle(0.45);
    \draw[color=gray,very thin,densely dashed] (0.6,0.6) circle (0.45);
    \draw[-, color=blue!60!black, very thick] (0.6,0.6) circle (0.03);
    \draw[color=blue!50!darkgray, ->, thick, densely dashed] (0.6,0.6) -- (0.93,0.93);
    \node[color=blue!60!black] at (0.72,0.90) {\footnotesize $r$};
    \node at (0.60,0.3) {{\footnotesize gas}};
    \node at (0.75,3.4) {Vessel};
    \node at (1.9,3.0) {Fiber};
    \draw[->, color=darkgray, very thick] (2.5,2.5) -- (2.5,4.5) node[color=black,pos=0.9,right] {$\pmb{y}$}; 
    \draw[-, color=blue!60!black, very thick] (2.5,2.5) circle (0.03);
    \coordinate (SS) at (2.7,1.4);
    \coordinate (EE) at (3.4,1.1);
    \draw[color=blue, ->, line width=0.08cm, opacity=0.8] (SS) to [out=290,in=180,distance=0.5cm] (EE);
    \node[color=blue] at (3.6,1.1) {\footnotesize $\pmb{U}$};
    \draw[color=black, fill=cyan!30!white] (5.3,2.2) rectangle ++(0.5,0.5);
    \draw[color=black, densely dashed, fill=white] (5.55,2.45) circle (0.18);
    \node[color=black] at (6.1,2.45) {$\Yfast$};
    \draw[color=black, densely dashed] (5.3,1.5) rectangle ++(0.5,0.5);
    \draw[color=black, fill=cyan!60!gray] (5.55,1.75) circle (0.20);
    \draw[color=white, fill=white] (5.55,1.75) circle (0.08);
    \draw[color=yellow, fill=yellow] (5.55,1.75) circle (0.05);
    \node[color=black] at (6.1,1.75) {$\Yslow$};
    \node[color=black] at (5.2,0.20) {$\Ymicro$};
    \draw[color=blue, ->, line width=0.05cm, opacity=0.8] (4.50,0.73) -- (4.50, 1.00);
    \draw[color=blue, ->, line width=0.05cm, opacity=0.8] (4.42,0.65) -- (4.15, 0.65);
    \draw[color=blue, ->, line width=0.05cm, opacity=0.8] (4.58,0.65) -- (4.85, 0.65);
    \draw[color=blue, ->, line width=0.05cm, opacity=0.8] (4.50,0.57) -- (4.50, 0.30);
    \draw[color=blue, -, line width=0.05cm] (4.5,0.65) circle(0.07);
    \node[color=blue] at (4.25,0.25) {\footnotesize $\pmb{\Uroot}$};
  \end{tikzpicture}
  }{\includegraphics[width=0.40\linewidth]{tikzfigure4a}}%
  \ifthenelse{\boolean{@IsTikzPlotsOnly}}{%
  \begin{tikzpicture}
    \draw[color=white] (0,0) rectangle (5,5); 
    \node at (2.6,5.3) {(b) {\rmfamily\normalfont Stem cross-section
        $\Xmacro$ tiled with reference cells.}}; 
    \foreach \x in {0,1,...,9}
      \foreach \y in {0,1,...,9}{
        \draw[color=white,dotted,fill=cyan!40!white]          (0.5*\x,0.5*\y) rectangle ++(0.5,0.5); 
        \draw[color=cyan!40!white,fill=cyan!60!gray!80!white] (0.5*\x+0.25,0.5*\y+0.25) circle(0.215);
        \draw[color=brown!80!black, densely dotted]           (0.5*\x+0.25,0.5*\y+0.25) circle(0.13);
        \draw[color=cyan!60!gray,fill=white]                  (0.5*\x+0.25,0.5*\y+0.25) circle(0.1);
        \draw[color=white,fill=yellow!80!white]               (0.5*\x+0.25,0.5*\y+0.25) circle(0.07);
        \draw[color=gray!30!white,fill=yellow!80!white]       (0.5*\x+0.07,0.5*\y+0.07) circle(0.06);
    }    
    \begin{scope}[even odd rule]
      \clip (2.5,2.5) circle(2.5) (0,0) rectangle (5,5);
      \fill[white] (0,0) rectangle (5,5);
    \end{scope}
    \draw[color=white] (0,0) rectangle (5,5); 
    \draw[color=brown, ultra thick] (2.5,2.5) circle(2.5);
    \draw[color=brown, ultra thick, fill=white] (2.5,2.5) circle (1.0);
    \draw[->, thick, densely dashed, color=blue!60!black] (2.5,2.5) -- (4.27,4.27) node[pos=0.9,left] {$\Rtree$};
    \draw[->, thick, densely dashed, color=blue!60!black] (2.5,2.5) -- (3.45,2.8) node[pos=0.6,below] {$\Rsap$};
    \draw[->, color=darkgray, very thick] (2.5,2.5) -- (2.5,4.9)  node[pos=0.95,right,color=black] {$\pmb{x}$};
    \draw[-, color=blue!60!black, very thick] (2.5,2.5) circle (0.03);
    \node[color=black] at (4.5,0.6) {$\Xmacro$};
  \end{tikzpicture}
  }{\includegraphics[width=0.51\linewidth]{tikzfigure4b}}}
  \caption{(a) The reference cell $\Ymicro$ containing a fiber located
    at the center (the brown dashed line is the fiber wall) outside of
    which lies the vessel. For the purposes of the periodic
    homogenization process, an artificial boundary $\Gamma$ (outer
    dotted circle) is introduced that separates $\Ymicro$ into two
    sub-regions: $\Yslow$, a fiber--vessel overlap region where heat
    diffusion is slow (shaded in medium blue); and $\Yfast$, the outer
    portion of the vessel region where diffusion is relatively fast
    (light blue). A gas bubble of radius $r$ is depicted in the
    lower-left corner, and the root water source $U_r$ in the
    lower-right. (b) An annular sapwood cross-section is tiled
    periodically with copies of the reference cell from (a). Mature
    trees contain a non-conducting heartwood region extending out to
    radius $x=\Rsap$, whereas younger saplings may have no heartwood
    ($\Rsap=0$).}
  \label{fig:macro2}
\end{figure}

To complete the model description we must derive two additional
equations describing heat transport: the first capturing microscale
effects that arise from phase change within the fibers and vessels; and
the second governing macroscale effects throughout the xylem that are
driven by ambient temperature variations. However, the temperatures on
the two scales are tightly coupled and so that extra care must be taken
to properly account for the transfer of heat energy between the micro-
and macro-scales.  For this purpose, we apply the method of periodic
homogenization or two-scale convergence \citep{allaire-1992} which
posits that in a material such as sapwood having a clear separation of
scales, the detailed microscale dynamics can be represented by a simpler
problem defined on a reference cell $\Ymicro$. For reasons of
simplicity, $\Ymicro$ is typically assumed to have radial symmetry,
which is not the case for the microscale geometry in
Figure~\ref{fig:cell-geometry}c; nevertheless, we can still define a
modified reference cell that has the requisite symmetry.  Because the
essential freeze--thaw processes that govern pressure generation occur
within fibers, we choose a fiber-centric coordinate system with
radial variable $y$, in which the fiber is placed in the middle of a
square reference cell having side length $\varepsilon$ (see
Figure~\ref{fig:macro2}a). Since the vessel is so much larger than the
fiber, it appears simply as a sap reservoir from the fiber perspective.
Consequently, we consider the vessel as the region of the
reference cell lying outside the fiber, as depicted in
Figure~\ref{fig:macro2}a (which in some sense ``turns the vessel
inside-out'').  In order that this modified reference cell remains
consistent with the original fiber--vessel geometry in
Figure~\ref{fig:cell-geometry}c, we impose a simple volume constraint
\begin{linenomath}\begin{gather}
  \pi \big({R^v}\big)^2 + \pi \big({R^f}\big)^2 N = \varepsilon^2, 
  \label{eq:eps-geometry}
\end{gather}\end{linenomath}
which ensures that each reference cell captures the net influence of one
vessel and $N$ adjacent fibers.  Although Figure~\ref{fig:macro2}a
depicts the vessel gas as a circular region in the lower-left corner,
the gas is not strictly assigned to any physical location otherwise the
radial symmetry would be broken.  Instead, it is represented in terms of
the radius $r$ of the equivalent gas bubble along with the fraction of
gas in dissolved form. As a result, the equations from the previous
section remain identical despite this apparent change in reference cell
geometry.

Next we exploit the regular, quasi-periodic microstructure of sapwood
and view the stem as being constructed of a periodic array of reference
cells as pictured in Figure~\ref{fig:macro2}b. The macroscopic domain
$\Xmacro$ is the 2D stem cross-section, which is an annular cylinder
having outer radius $x=\Rtree$ and inner radius $x=\Rsap$ (bounding the
non-conductive heartwood).
Our aim is now to derive two equations for heat transport: one on the
reference cell that incorporates local variations in temperature due to
freeze--thaw processes occuring in a given reference cell at location
$x\in\Xmacro$; and a second equation capturing macroscopic heat
transport throughout $\Xmacro$ in response to ambient temperature
fluctuations, combined with the local effects.  The main objective of
the periodic homogenization process is to derive appropriate heat
transport coefficients for the macroscale heat equation that incorporate
the effects of the microscale by averaging the solution of the reference
cell problem $\Ymicro$ appropriately. Note that the size $\varepsilon$
of the reference cell plays a dual role in homogenization:
it can be considered as a physical dimension, but it must also be viewed
asymptotically in the limit as $\varepsilon\to 0$ to obtain the averaged
effect of the microscale freeze--thaw process on the macroscale.

We now summarize the essential aspects of the periodic homogenization
procedure, for which complete mathematical details are provided in
\citet{konrad-peter-stockie-2017}. For technical reasons, the reference
cell $\Ymicro$ is separated into two sub-regions ($\Yslow$ and $\Yfast$)
pictured in Figure~\ref{fig:macro2}a, where $\Yslow$ refers to the fiber
plus the inner portion of the vessel where heat diffuses slowly, whereas
$\Yfast$ is the remaining outer portion of the vessel where diffusion is
fast compared to $\Yslow$. The curve $\Gamma$ is an artificial boundary
separating $\Yfast$ from $\Yslow$ so that
$\Ymicro(x,t)=\Yfast\,\cup\,\Gamma\,\cup\,\Yslow(x,t)$.  Note that an
implicit time and space dependence appears in $\Yslow$ (and hence also
$\Ymicro$) owing to the motion of phase boundaries that
alters the geometry of $\Ymicro$ depending on the specific location $x$
(although our notation will often omit this dependence).  Our aim is
then to derive two heat diffusion equations, one for $\Tmicro(x,y,t)$ on
the microscale domain $\Yslow \times \Xmacro$ and the other for
$\Tmacro(y,t)$ on the macroscale domain $\Xmacro$.

The governing equations are stated in a mixed temperature--enthalpy
formulation in order to properly capture phase transitions.  To this end
we define $\Tmacro(x,t)$ and $\Emacro(x,t)$ as the macroscale
temperature and enthalpy variables, which are both constant inside
$\Yfast$ and thus depend only on the macroscale spatial coordinate $x$.
The corresponding microscale variables on $\Yslow(x,t)$ are
$\Tmicro(x,y,t)$ and $\Emicro(x,y,t)$, which vary at each point in the
macroscopic domain as well as the microscale $y$.  We impose the usual
temperature--enthalpy relationship on both $\Tmacro=\omega(\Emacro)$ and
$\Tmicro=\omega(\Emicro)$ where
\begin{linenomath}\begin{gather}
  \omega(\Emacro) = 
  \begin{cases}
    \frac{\Emacro}{c_i}, & \quad \text{if $\Emacro < E_i-\delta_i$},\\
    \Tcrit-\frac{2\Emacro-E_i-E_w}{2c_\infty},  
    & \quad \text{if $E_i-\delta_i \leqslant \Emacro \leqslant E_w+\delta_w$},\\ 
    \Tcrit+\frac{\Emacro-E_w}{c_w}, 
    & \quad \text{if $E_w+\delta_w< \Emacro$},
  \end{cases}
  \label{eq:TE-omega} 
\end{gather}\end{linenomath}
captures implicitly the change in phase that occurs at the
\mymod{critical (freezing)} temperature $\Tcrit$.  This $\omega$ is a
piecewise function consisting of two linear segments with slopes
$c_i^{-1}$ in ice and $c_w^{-1}$ in liquid, connected by a steep layer
with slope $c_\infty^{-1}$ (where $c_\infty \approx 10^7$).  The
constants $c_w$ and $c_i$ refer to the specific heat capacities of water
and ice, while $\delta_i=\frac{c_i(E_w-E_i)}{2(c_\infty-c_i)}$ and
$\delta_w=\frac{c_w(E_w-E_i)}{2(c_\infty-c_w)}$ are chosen to ensure
$\omega(E)$ is continuous.

Heat transport within the liquid-filled subregion $\Yslow$ is governed by
the heat equation 
\begin{linenomath}\begin{gather}
  c_w \partial_t \Tmicro - \yderiv{ \Big(D(\Emicro) \yderiv{\Tmicro}
    \Big)} = 0  \qquad \text{in $\Yslow(x,t) \times \Xmacro$},
  \label{eq:temp-micro}
\end{gather}\end{linenomath}
where $D$ is a thermal diffusion coefficient that is also
a given piecewise linear function of enthalpy \citep{visintin-1996}:
\begin{linenomath}\begin{gather}
  D(\Emicro) =
  \begin{cases}
    \frac{k_i}{\myrho_i}, & \quad \text{if $\Emicro < E_i$},\\
    \frac{k_i}{\myrho_i}+\frac{\Emicro-E_i}{E_w-E_i} \Big(\frac{k_w}{\myrho_w}
    - \frac{k_i}{\myrho_i}\Big), 
    & \quad \text{if $E_i \leqslant \Emicro \leqslant E_w$},\\ 
    \frac{k_w}{\myrho_w}, & \quad \text{if $E_w < \Emicro$}.
  \end{cases}
  \label{eq:diffusivity}
\end{gather}\end{linenomath}
The partial differential equation \eqref{eq:temp-micro} requires
boundary conditions on the inner ($\partial\Yslow$) and
outer ($\Gamma$) boundaries of $\Yslow$, which are
\begin{linenomath}\begin{align}
  \Tmicro &= \Tcrit  \qquad \text{on $\partial \Yslow(x,t) \times
    \Xmacro$ \quad (phase-change boundary)},
  \label{eq:microBC-inner}\\
  \Tmicro &= \Tmacro \qquad \text{on $\Gamma \times \Xmacro$ 
    \quad (coupling to macro-temperature)}. 
  \label{eq:microBC-outer}
\end{align}\end{linenomath}

The homogenization procedure that we apply next to obtain a macroscale
heat equation is more complicated and requires first taking
$y=x/\varepsilon$ in the microscale problem and then expanding the
solution asymptotically as $\varepsilon\to 0$, which is referred to as
the \emph{two-scale limiting process}. The governing equation is written
in an integral (weak) formulation, but after approximating the resulting
integral terms in the $\varepsilon\to 0$ limit one obtains the following
strong formulation of the limit problem:
\begin{linenomath}\begin{gather}
 \partial_t \Emacro - 
  \partial_x \Big(\Pi D(\Emacro)\partial_x \Tmacro\Big)
  = \frac{1}{|\Yfast|} \int_{\Gamma} D(\Emicro) \nderiv{\Tmicro} \, dS
  \qquad \text{in $\Xmacro$} .
  \label{eq:temp-macro} 
\end{gather}\end{linenomath}
Note that this is an alternate form of the heat diffusion equation,
written in a mixed temperature--enthalpy form that implicitly captures
parameter discontinuities across phase interfaces.  There are two new
terms appearing in Eq.~\eqref{eq:temp-macro} via the homogenization
process that are critically important in properly capturing the
influence of the microscale problems on the macroscale:
\begin{itemize}
\item The homogenized diffusion operator contains an extra constant
  factor $\Pi$, which is a purely geometric quantity consisting of a
  $2\times 2$ matrix with entries
  
  \begin{linenomath}\begin{gather}
    \Pi_{ij} = \frac{1}{|\Yfast|}
    \int_{\Yfast} (\delta_{ij} + \yderiv{\mu_i}) \, dy.
    \label{eq:Pi-matrix}
  \end{gather}\end{linenomath}
  Here, $\delta_{ij}$ is the Kronecker delta symbol and $\mu_i(y)$ are 
  solutions to simple elliptic PDE problems on the reference cell
  $\Yfast$ \citep{allaire-1992}.
\item The source term on the right-hand side is a surface integral over
  the artificial boundary $\Gamma$ of the microscopic heat
  flux. 
\end{itemize}
It is important to recognize here that the temperature--enthalpy
relationship $T=\omega(E)$ in \eqref{eq:TE-omega} involves the
critical temperature, $\Tcrit$.  Within the fiber (which contains pure
water) we take $\Tcrit=0\degC$, but in the macroscale problem $\Tcrit$
must be replaced with
\begin{linenomath}\begin{gather}
  \TcSap = \Tcrit - \frac{K_b C_s}{\myrho_w},
  \label{eq:FPD}
\end{gather}\end{linenomath}
where $K_b=1.853$ is the cryoscopic (or Blagden's) constant.  This
accounts for the freezing point depression or FPD that arises due to
dissolved solutes (primarily sugar) in the vessel sap.  Finally, the PDE
\eqref{eq:temp-macro} is supplemented with the boundary condition
\begin{linenomath}\begin{gather}
  \Tmacro = \Tamb(t)  \qquad \text{on $\partial \Xmacro$},
  \label{eq:macroBC}
\end{gather}\end{linenomath}
which sets the outer stem surface temperature equal to a given ambient
air temperature and is what ultimately drives the freeze--thaw
process.  Complete details of the homogenization procedure can be found
in \citet{konrad-peter-stockie-2017}, and we also refer the interested
reader to the work of \citet{chavarriakrauser-ptashnyk-2013} who applied
periodic homogenization to a related problem involving osmotic transport
in non-woody plants.

The parameter values appearing in these equations are taken mostly from
previous work \citep{graf-ceseri-stockie-2015,
  konrad-peter-stockie-2017} that focused on comparisons to experimental
data from black walnut. These parameters are summarized in
Table~\ref{tab:params}, with a few small adjustments for
red/sugar maple as indicated.  The ``base case'' that is indicated there
corresponds to a sugar maple sapling with diameter of 14~cm that has sap
sugar content of 3\% by mass.


\begin{table}[btp]
  \footnotesize
  \centering
  \caption{Model parameters used in the base case simulations are taken from
    \citet{graf-ceseri-stockie-2015}, unless otherwise indicated.
    Modifications to parameters $\Rtree$, $\Rsap$, $\gamma_s$ for
    comparison with the red/sugar maple experiments are detailed in the
    text.}  
  \label{tab:params}
  \newcommand{\mytwocol}[2]{\parbox{\widthof{0.070,0.093}}{\raggedleft
      #1}~~{\color{gray!50}$\mid$}~~\parbox{\widthof{0.070,0.093}}{\raggedright #2}}
  \setstretch{1.04} 
  \begin{tabular}{clcc}\hline
    {\bf Symbol} & {\bf Description} & {\bf Values} & {\bf Units} \\\hline
    & & & \\[-0.2cm]
    \multicolumn{4}{l}{\emph{Variables (functions of time $t$ and space
        $x$ or $y$):}}\\
    $s_{iw}$, $s_g$ & fiber interface locations & & m \\ 
    $r$      & vessel bubble radius & & m \\
    $U$      & water transferred from fiber to vessel & & m${}^3$ \\
    $\Uroot$ & root water uptake & & m${}^3$ \\
    $\Tmicro$, $\Tmacro$ & temperature & & \degC \\
    $p$      & pressure (relative to atmospheric) & & Pa \\
    $\myrho$   & density  & & kg m${}^{-3}$ \\
    $V$      & volume   & & m${}^3$ \\
    & & & \\[-0.2cm]
    \multicolumn{4}{l}{\hspace*{0.5cm}\emph{(Subscripts:  $i$,$w$,$g$ for ice,
        water/sap, gas; Superscripts: $f$,$v$ for fiber, vessel)}}\\
    & & & \\[-0.2cm]
    \multicolumn{4}{l}{\emph{Tree physiological parameters:}}\\
    $\Afvwall$ & surface area of fiber--vessel wall & \myee{6.28}{-8} & m${}^2$ \\
    $\Aroot$    & root area per vessel $=\Atree(R^v/\Rtree)^2$ & \myee{1.14}{-6} & m${}^2$ \\
    $\varepsilon$ & side length of reference cell, Eq.~\eqref{eq:eps-geometry} & 
    \myee{4.33}{-5} & m \\
    $L^f$       & length of fiber & \myee{1.0}{-3} & m \\
    $L^v$       & length of vessel element & \myee{5.0}{-4} & m \\
    $\tyreecond$& hydraulic conductivity of fiber--vessel wall & \myee{5.54}{-13} & m\,s${}^{-1}$\,Pa${}^{-1}$ \\
    $\tyreecond_r$ & root hydraulic conductivity
                     \mymod{\citep{steudle-peterson-1998}} &
    \myee{2.7}{-16} & m\,s${}^{-1}$\,Pa${}^{-1}$ \\
    %
    $N$         & number of fibers per vessel & $16$ & -- \\
    $R^f$       & inside radius of fiber & \myee{3.5}{-6} & m \\
    $R^v$       & inside radius of vessel & \myee{2.0}{-5} & m \\
    & & & \\[-0.2cm]
    \multicolumn{4}{l}{\emph{Physical constants:}}\\
    $\Henry$    & Henry's constant for air dissolved in water & $0.0274$ & -- \\
    $K_b$       & cryoscopic (Blagden) constant & $1.853$ & kg\,\degC\,mol${}^{-1}$ \\
    $M_g$       & molar mass of gas (air) & $0.029$ & kg\,mol${}^{-1}$ \\
    $M_s$       & molar mass of sugar (sucrose) & $0.3423$ & kg\,mol${}^{-1}$ \\
    $\Rgas$     & universal gas constant & $8.314$ & J\,\degC${}^{-1}$\,mol${}^{-1}$ \\
    & & & \\[-0.2cm]
    \multicolumn{2}{l}{\emph{Water phase properties:}} & 
    \mytwocol{\emph{Ice}}{\emph{Water}} & \\
    $c_i$, $c_w$ & specific heat capacity & \mytwocol{$2100$}{$4180$}
    & J\,\degC${}^{-1}$\,kg${}^{-1}$ \\
    $E_i$, $E_w$ & enthalpy at $\Tcrit$ & \mytwocol{$574$}{$907$} & kJ\,kg${}^{-1}$ \\
    $k_i$, $k_w$ & thermal conductivity & \mytwocol{$2.22$}{$0.556$}
    & W\,m${}^{-1}$\,\degC${}^{-1}$ \\
    %
    %
    %
    $\myrho_i$, $\myrho_w$ & density & \mytwocol{$917$}{$1000$} & kg\,m${}^{-3}$ \\
    $\sigma_{iw}$, $\sigma_{gw}$ & surface
    tension \citep{fowler-krantz-1994} & \mytwocol{$0.033$}{$0.076$} & N\,m${}^{-1}$ \\
    $c_\infty$   & enthalpy regularization parameter, Eq.~\eqref{eq:TE-omega} & 
    \myee{1.0}{7} & J\,\degC${}^{-1}$\,kg${}^{-1}$ \\
    & & & \\[-0.2cm]
    \multicolumn{2}{l}{\emph{Base case parameters:}} & & \\
    $\Rtree$    & tree stem radius & $0.07$ & m \\
    $\Rsap$     & sapwood/heartwood boundary & $0$ & m \\
    $\heartfrac$& heartwood fraction $=\Rsap/\Rtree$ & 0 & --\\
    $\gamma_s$  & sap sugar content (mass fraction) & $0.03$ & -- \\
    $\Atree$    & total root area \citep{day-harris-2007} & 14 & m${}^2$ \\
    $C_s$       & sap sugar concentration $= \gamma_s\myrho_w/M_s$ & $87.6$ & mol\,m${}^{-3}$ \\
    $\TcSap$    & freezing point depression (FPD) $= - K_b C_s/\myrho_w$ & $-0.162$ & \degC \\
    \mymod{$\Psi_s$} & \mymod{soil water potential $=p_w^v(0)$} & \mymod{\myee{2.03}{5}} & \mymod{Pa} \\
    \hline
  \end{tabular}
\end{table}

\subsection{Numerical solution algorithm}
\label{sec:numerics}

The exudation model equations consist of five ODEs
\eqref{eq:sgas}--\eqref{eq:Uroot-darcy} for the microscale variables
along with two PDEs \eqref{eq:temp-micro} and \eqref{eq:temp-macro} for
temperature.  We apply the method-of-lines to discretize the temperature equations
by first approximating all spatial derivatives using a finite volume 
approach, which yields a large coupled system of time-dependent
ODEs. The macroscale variable $x$ is discretized at $n_x$ equally-spaced
points, with $n_x$ chosen between 25--50 (depending on tree size) so
that the grid spacing $\Delta x=(\Rtree-\Rsap)/n_x$ is less than 0.3~cm,
which we find is sufficient in practice to resolve the freezing and thawing
fronts. For the microscale problem, we obtain satisfactory accuracy with
a relatively coarse grid having $n_y=6$ points.

Assembling the semi-discrete equations for temperature together with the
remaining ODEs and algebraic constraints within each reference cell at
location $x$ yields a coupled differential--algebraic system that is
integrated in time using the stiff solver {\tt ode15s}
in~\citet{matlab-2020a}.  We note that a stiff solver is required for
this problem because of the widely disparate time scales arising from
the disparate dynamics on the cellular level and within the
tree stem.  The strong coupling between micro- and macroscale
temperatures is handled by applying a split-step time discretization:
first, the microscale equations for $s_g$, $s_{iw}$, $r$, $U$, $\Uroot$,
$\Tmicro$ in the reference cell at each discrete point are advanced
in time while holding the macroscale temperature $\Tmacro$ constant;
following that, $\Tmacro$ is advanced in time by holding other
variables constant. More details related to implementation of the
multiscale algorithm can be found in \citet{graf-ceseri-stockie-2015}.

\section{Results}
\label{sec:results}

\subsection{Experimental data and temperature smoothing}
\label{sec:measurements}

\begin{figure}[tbhp]
  \centering\small
  \ifthenelse{\boolean{@IsProofs}}%
  {\includegraphics[width=0.9\textwidth]{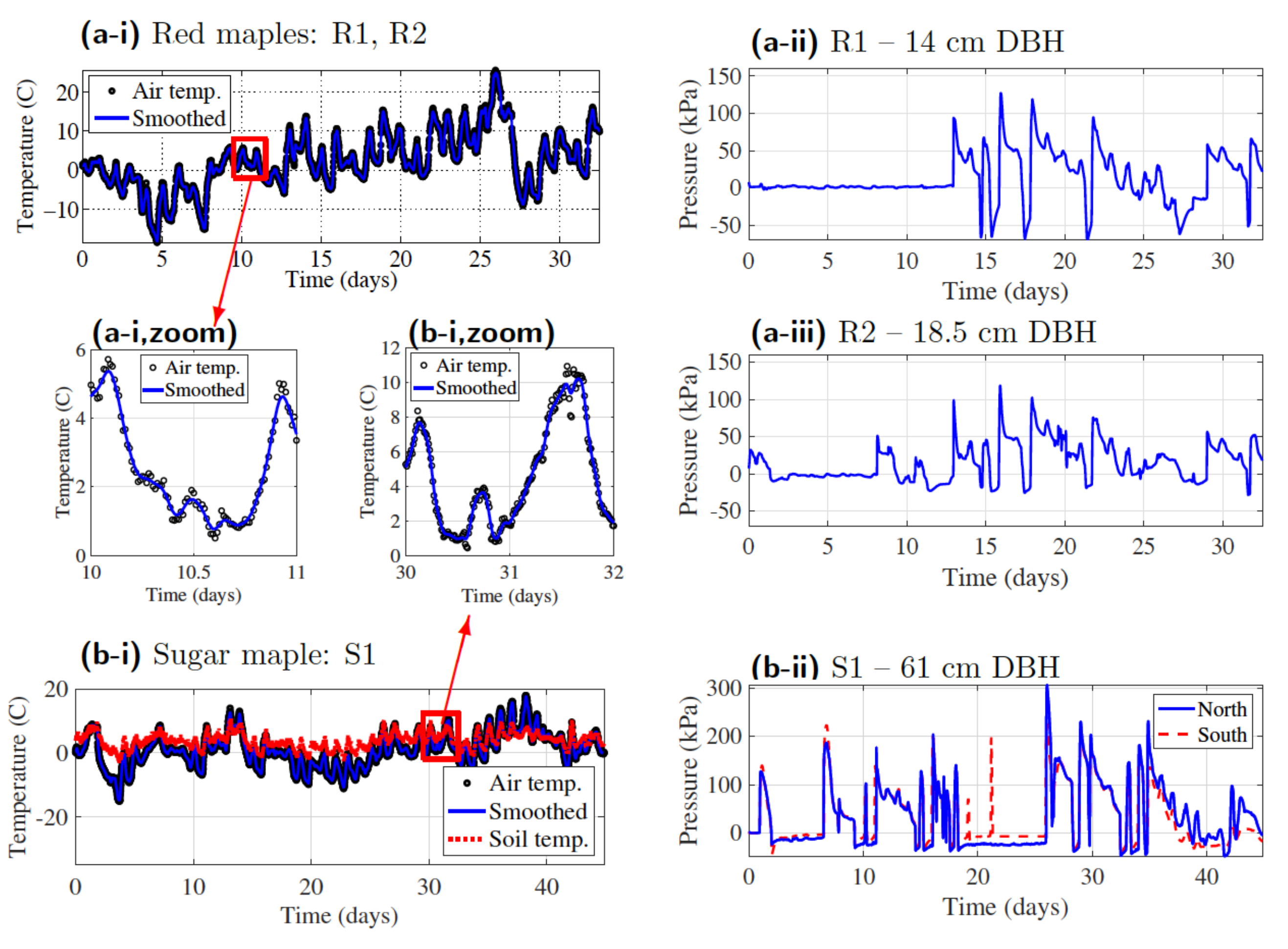}}%
  {%
  \begin{tabular}{lcl}
    && \hspace*{0.8cm}\raisebox{0.05cm}{\myfiglett{(a-ii)} R1 -- 14~cm DBH} \\ 
    \multirow{3}{*}[2.8cm]{%
      {\setlength{\unitlength}{0.4\textwidth}
        \begin{picture}(1,1.5)
          \put(0,1.48){\hspace*{0.8cm}\myfiglett{(a-i)} Red maples: R1, R2}
          \put(0,1.05){\includegraphics[width=0.41\linewidth]{uvm-data/redtempsmzoom}}
          \put(0.08,0.52){\includegraphics[width=0.17\linewidth]{uvm-data/redtempsmzoom2}}
          \put(0.60,0.52){\includegraphics[width=0.18\linewidth]{uvm-data/sugartempsmzoom2}}
          \put(0,0.42){\hspace*{0.8cm}\myfiglett{(b-i)} Sugar maple: S1}
          \put(0,-0.03){\includegraphics[width=0.41\textwidth]{uvm-data/sugartempsmzoom}}
          \put(0.14,0.97){\myfiglett{(a-i,zoom)}}
          \put(0.68,0.97){\myfiglett{(b-i,zoom)}}
          {\linethickness{2mm}\thicklines\color{red}
            \put(0.415,1.25){\vector(-1,-4){0.063}}
            \put(0.74,0.32){\vector(1,4){0.046}}}
        \end{picture}
      }}
    && \includegraphics[width=0.4\textwidth]{uvm-data/redpresR1}\\
    && \hspace*{0.8cm}\myfiglett{(a-iii)} R2 -- 18.5~cm DBH\\
    && \includegraphics[width=0.4\textwidth]{uvm-data/redpresR2}\\[0.55cm]
    && \hspace*{0.8cm}\myfiglett{(b-ii)} S1 -- 61~cm DBH\\[-0.1cm]
    && \includegraphics[width=0.4\textwidth]{uvm-data/sugarpresS1}\\
  \end{tabular}}
  \caption{(Left, a-i and b-i) Measured air temperatures are plotted for
    two red maple trees (R1, R2) and one sugar maple (S1) from the UVM
    experiments. The temperature plot for sugar maple in (b-i) also
    includes values of soil temperature (at 30~cm depth) which remain
    mostly above 0\degC, hence supporting the assumption that liquid
    water is present even when air temperatures are below freezing.  The
    raw temperature data (blue points) are regularized by applying a
    simple weighted-average smoothing -- the resulting smoothed data are
    shown alongside the original temperatures in the zoomed views
    (a-i,zoom, b-i,zoom).  (Right, a-ii, a-iii and b-ii) Corresponding
    pressure data for the three trees.  An extra set of pressure
    measurements is included in the sugar maple plot (b-ii) to
    illustrate the impact of taking measurements on the north/south
    sides of the stem.}
  \label{fig:tree-data}
\end{figure}

Based on the experiments described in Section~\ref{sec:experiments}, we
have chosen to focus our attention on two red maple trees (which we
label R1 and R2 for convenience) and one sugar maple (labelled S1).  A
more detailed description of these experiments that emphasizes the
temperature measurements can be found in \citet{wilmot-2006}.  Air
temperature data for all three trees are displayed in
Figure~\ref{fig:tree-data}(left) over periods of 33 days (R1,R2) and 45
days (S1). These samples were singled out for comparison with numerical
simulations because the air temperature in each case features several
pronounced freeze--thaw cycles during the measurement period.  Note that
the temperature plot for sugar maple in Figure~\ref{fig:tree-data}b-i
includes a second curve (red, dashed) showing the soil temperature at
30~cm depth which for most of the 45-day period remains positive, even
during times when the air temperature is below zero.  This provides
strong evidence in support of the earlier
assumption~\ref{assume:soil-liquid} that liquid water is available for
root uptake or water extrusion even under freezing conditions.
Water can be extruded from roots to soil if the sum of
hydrostatic pressure ($\myrho g h$) and stem pressure exceeds the
osmotic pressure of sugar and solutes in the xylem of minor roots.

Although these air temperature measurements have an inherent large-scale
oscillation that varies roughly on a daily period, the two zoomed-in
views in Figures~\ref{fig:tree-data}a,b-i show that there are also
significant fluctuations from one 15-minute time interval to the next.
These rapid changes are likely due to a combination of local temperature
variability and measurement errors and are a major distinguishing
feature that sets these field measurements apart from others obtained
under carefully controlled laboratory conditions.  Indeed, just such an
experiment on black walnut trees \citep{ameglio-etal-2001} was used to
validate a previous incarnation of our multiscale model
\citep{graf-ceseri-stockie-2015} where the input temperatures was
specified as a given smoothly varying function of time.  Consequently,
our measurements for red/sugar maple provide an excellent opportunity to
validate the model under more realistic conditions.

Because our exudation model is based on differential equations that
expect the ambient temperature $\Tamb(t)$ in Eq.~\eqref{eq:macroBC} to
vary continuously in time, we need to impose some regularization to
smooth the raw temperature data. To this end, we apply a simple
weighted-average smoothing procedure in which each temperature value is
averaged with its two neighboring points using weights
$\big[\frac{1}{4}, \frac{1}{2}, \frac{1}{4}\big]$, with this procedure
being repeated 10 times. The smoothed temperature is displayed as a
solid curve along with the raw data in the two zoomed plots in
Figures~\ref{fig:tree-data}a,b-i, from which it is clear that this
procedure eliminates many irregularities without sacrificing much
detail. There is of course a risk that genuine fine-scale variations in
temperature are suppressed, but we have observed that reducing the
number of smoothing steps has no appreciable effect on the model
simulations.  This is consistent with the results in
Section~\ref{sec:results-sims} which show that exudation behavior is
dominated by the location of temperature zero-crossings and influenced
much less by variations in $\Tamb(t)$ away from zero.

\subsection{Parameter sensitivity study}
\label{sec:results-sensitivity}

To study the relative importance of various geometric and physical
properties on the sap exudation process, we identify a base case using
the parameters listed in Table~\ref{tab:params} and then vary certain
parameters relative to these base values.  Our chosen base case
represents a young sapling with $\Rtree=0.07$~m and $\Rsap=0$ that
hasn't yet developed any heartwood, and we use our best estimates for
the remaining parameters. To mimic a repeated sequence of diurnal
freeze--thaw events, we impose a simple sinusoidally varying ambient
temperature, $\Tamb(t)=5-15\sin\left({2\pi t}/{86400}\right)$,
that oscillates between $-10$ and $+20$\degC\ over a time interval of 5
days. This is admittedly an extreme range of temperatures, but it
does ensure that the entire stem is able to freeze and thaw completely
during each cycle. Results of the parameter sensitivity study are
presented in Figure~\ref{fig:sensitivity} as plots of root water uptake
$\Uroot$ and stem-averaged pressure $\overline{p} = \frac{1}{|\mathcal{A}|}
\int_{\mathcal{A}} p_w^v\, d{A}$, where $|\mathcal{A}|$ is the area of
the annular-shaped sapwood region.  This averaged pressure is a better
representation than any specific point value for what is measured by a
pressure gauge inserted into a taphole.

\newcommand{\mytabrow}[3]{%
  \begin{minipage}[t]{0.14\textwidth}
    \begin{center}
      \mbox{}\break
      {\bfseries\sffamily (#1)} \break 
      #2
    \end{center}
  \end{minipage} & 
  \begin{minipage}[t]{0.38\textwidth}
    \vspace*{0pt}
    \includegraphics[width=\textwidth]{maryam-code/new_11Sep2021/baseruns_MEL/#3_pres_compare}
  \end{minipage} & &
  \begin{minipage}[t]{0.38\textwidth}
    \vspace*{0pt}
    \includegraphics[width=\textwidth]{maryam-code/new_11Sep2021/baseruns_MEL/#3_uroot_compare}
  \end{minipage}}

\begin{figure}[tbhp]
  \centering\footnotesize
  \ifthenelse{\boolean{@IsProofs}}%
  {\includegraphics[width=\textwidth]{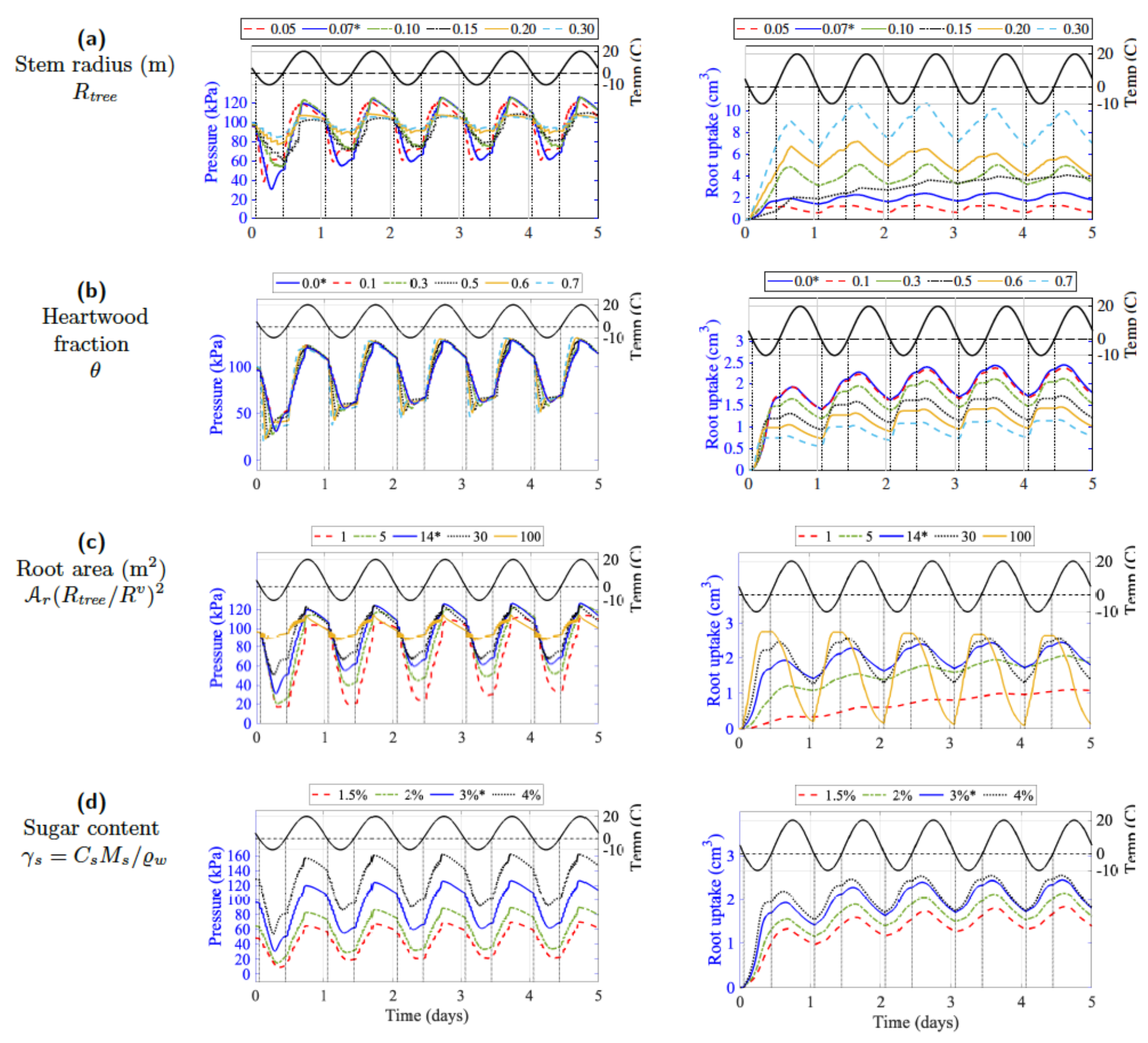}}%
  {%
  \begin{tabular}{cccc}
    \mytabrow{a}{Stem radius (m)\break $\Rtree$}{rtree2} \\
    & & & \\[-0.2cm]
    \mytabrow{b}{Heartwood fraction \break $\heartfrac$}{heartfrac2} \\
    & & & \\[-0.2cm]
    \mytabrow{c}{Root area (m${}^2$) \break $\Aroot (\Rtree/R^v)^2$}{aroot2} \\
    & & & \\[-0.2cm]
    \mytabrow{d}{Sugar content \break $\gamma_s=C_s M_s/\myrho_w$}{sugar2}
  \end{tabular}}
  \caption{Parameter sensitivity study showing the effects on
    stem-averaged pressure ($\overline{p}$) and root water uptake
    ($\Uroot$) due to variations in four parameters:
    (a) stem radius,
    (b) heartwood fraction,
    (c) total root area,
    and (d) sap sugar content.  In each plot, the parameter value from
    the base case is highlighted in the legend with a ``$\ast$'' and the
    corresponding curve is drawn with a blue solid line.  The dotted
    vertical lines indicate freeze and thaw events, which are times when
    the ambient temperature crosses zero.  \mydel{The ``envelope
      curves'' tracing max-min points for the $\Rtree$ base-case
      pressure are shown as two dashed blue lines (a, left).}}
  \label{fig:sensitivity}
\end{figure}

The results plotted in Figure~\ref{fig:sensitivity} focus on
variations in four key model parameters and how they affect the
behavior of pressure and root uptake:
\begin{enumerate}[label=(\alph*)]
\item Stem radius ($\Rtree$), which is the primary geometric parameter
  that distinguishes between the mature trees in this study and younger
  saplings.  Values of $\Rtree$ are selected between 5--30~cm which
  covers the range of tree sizes in the experiment discussed in
  Section~\ref{sec:experiments} and the corresponding $\overline{p}$ and
  $\Uroot$ solution curves are shown in Figure~\ref{fig:sensitivity}a.
  Note that the curve corresponding to the base case ($\Rtree=0.07$) is
  always drawn as a solid blue line and is highlighted in the legend
  with a ``$\ast$''.

\item Heartwood fraction ($\heartfrac=\Rsap/\Rtree$), which is zero for
  young saplings but can be significantly larger in mature trees.  The
  red/sugar maples considered in this study are from a well-established
  area of the forest in which trees typically have 25--50\%\ of their
  basal area taken up by heartwood (with lower fractions in smaller
  trees and higher fractions in larger trees).  Most of these trees have
  been tapped annually for maple collection over a period of 50--60
  years, which generates a column of non-conductive wood that extend
  above and below each year's taphole (typically by a distance of
  $\pm$0.25\,m). Therefore, while there is undoubtedly some heartwood in
  these trees, there is also considerable non-conductive wood within the
  tapping band as a result of tapping history.  We therefore chose
  values of heartwood fraction within the range $\heartfrac\in[0,0.7]$,
  which is consistent with the measurements of
  \citet{duchesne-etal-2016} who found a maximum value of
  $\heartfrac\approx 0.45$, while also allowing for even higher values
  such as those reported by \citet{baral-etal-2017}.

\item Total root area ($\Atree=\Aroot(\Rtree/R^v)^2$), which is related
  to the root area per vessel ($\Aroot$) by scaling proportionally to
  cross-sectional area.  Because values of $\Atree$ for red maple have
  been reported to lie in the range 10.4 to 18.6~m${}^2$
  \citep{day-harris-2007}, we choose a value of $\Atree=14$ for the base
  case that lies near the middle of this range; scaling by the area
  ratio yields a corresponding root area per vessel of $\Aroot=
  \Atree(R^v/\Rtree)^2 \approx 1.14\times 10^{-6}$~m${}^2$. For the
  sensitivity results shown in Figure~\ref{fig:sensitivity}c, we have
  actually selected a wider range of $\Atree\in[1, 100]$.

  \leavethisout{
  \item Root reflection coefficient for outflow ($\Cout$), which is the
    major extension we have made to the original exudation model
    that assumed $\Cout\equiv 0$ \citep{graf-ceseri-stockie-2015}.  There
    is an extensive literature suggesting that root conductivity in a wide
    range of trees and plants is not constant but rather exhibits both
    seasonal and diurnal variations controlled by aquaporin
    membranes within the roots~\citep{javot-maurel-2002, steudle-1994}.
    Furthermore many root systems exhibit an asymmetry in conductivity
    between inflow and outflow, which can be modelled by means of a
    reflection coefficient $\Cout$~\citep{knipfer-fricke-2010,
      steudle-1994} that is less than 1 and can drop to $\Cout=0.2$ or
    less~\citep{henzler-etal-1999}. Among the few studies we are aware of
    that mention maple trees, \citet{oleary-1965} observed no appreciable
    root outflux ($\Cout\approx 0$) whereas \citet{dawson-1993} observed
    more moderate asymmetry (with $\Cout$ significantly greater than 0).
    In the absence of any reliable estimates specific to maple, we assume
    a base value of $\Cout=0.2$ that is constant in time, and compare with
    other values from the range $[0.05, 0.75]$ in
    Figure~\ref{fig:sensitivity}d. }
  
\item Sap sugar content by mass ($\gamma_s$), which is related to sugar
  concentration by $\gamma_s=C_sM_s/\myrho_w$, where $M_s$ the molar mass
  of sugar.  High sugar content is an important feature distinguishing
  maples from other species that exude sap.  Sugar maple sap contains
  roughly 3\%\ sugar on average during the sap harvest season
  \citep{larochelle-etal-1998} but can be as high as 5\%\ in some trees
  \citep{jones-alli-1987}. On the other hand, red maples tend to have a
  lower sugar content that is closer to 2\%\ on average.  Sugar content
  also 
  varies significantly between seasons, between trees, and also
  throughout a given season (starting a bit low, rising for the first
  1/4 to 1/3 of the season, then steadily dropping towards the end).  We
  have therefore chosen a representative value of 3\%\ ($\gamma_s=0.03$)
  for the base case along with several other values selected from the
  range 1.5 to 4.0\%\ as depicted in Figure~\ref{fig:sensitivity}e.
\end{enumerate}
This parameter sensitivity study is partly based on results from
\citet{zarrinderakht-mscthesis-2017}, which includes additional results
not reported here.


\subsection{Numerical simulations of red and sugar maple}
\label{sec:results-sims}

\begin{figure}[tbhp]
  \centering\small
  \ifthenelse{\boolean{@IsProofs}}%
  {\includegraphics[width=0.75\textwidth]{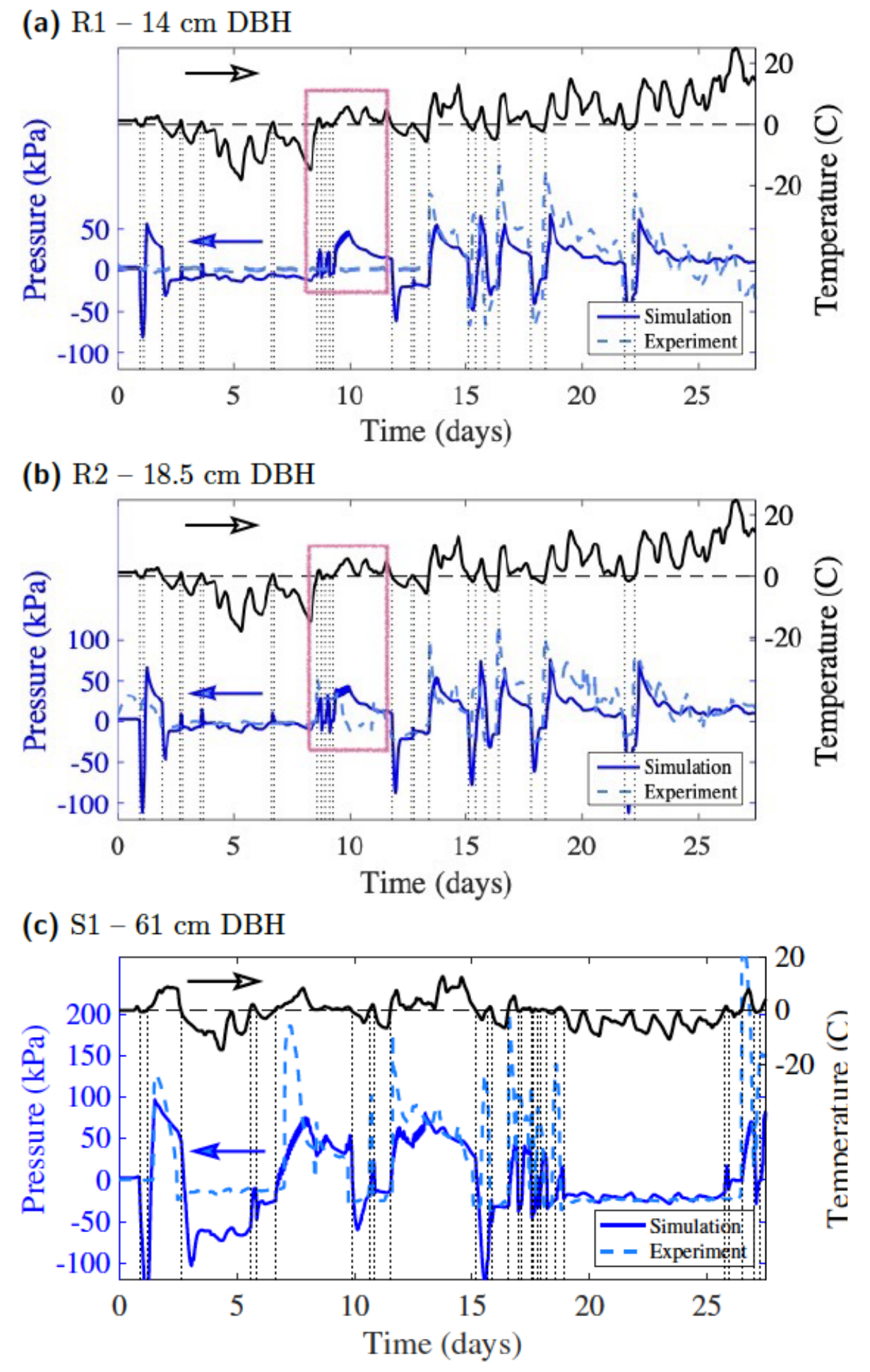}}%
  {%
  \begin{tabular}{ll}
    \myfiglett{(a)} R1 -- 14~cm DBH & \\
    \includegraphics[width=0.6\textwidth]{Mel_R1_prestempR1_mod}\\
    %
    \myfiglett{(b)} R2 -- 18.5~cm DBH & \\
    \includegraphics[width=0.6\textwidth]{Mel_R2_prestempR2_mod}\\
    %
    \myfiglett{(c)} S1 -- 61~cm DBH & \\
    \includegraphics[width=0.6\textwidth]{Mel_S1_prestempS1}
  \end{tabular}}
  \caption{Comparison of stem pressure from experiments (solid blue
    lines) and simulations (dashed cyan lines) for red maple (a,b) and
    sugar maple trees (c).  The smoothed temperature data are also
    displayed on the top set of axes, and each freeze--thaw event is
    highlighted with a vertical dotted line (at each time where
    temperature crosses 0\degC) so that these events are easily
    connected with corresponding pressure spikes.  A ``weak thaw'' event
    is highlighted with a pink box in the R1,R2 simulations (a,b).}
  \label{fig:sims}
\end{figure}

We next apply the MATLAB code to simulate the three trees
singled out in Section~\ref{sec:measurements}, taking the smoothed
temperature curves depicted in Figures~\ref{fig:tree-data}a,b-i as input
for the ambient temperature $\Tamb(t)$.  Two simulations are performed
for red maple trees (R1,R2) with stem radii $\Rtree=7$ and $9.25$~cm,
both consisting entirely of sapwood ($\Rsap=0$) since they are relatively
young trees.  The sugar content for these trees is set to 1.8\%
which is a representative value for mid-to-late season, but otherwise
all model parameters are the same as the base case in
Table~\ref{tab:params}.  The equations were integrated over a period of
27~days, which covers the majority of freeze--thaw events occurring in
the air temperature measurements. The resulting plots of simulated average
pressure $\overline{p}$ are displayed in Figures~\ref{fig:sims}a,b
alongside the corresponding experimental measurements. The smoothed
temperature data are shown at the top of each plot, with dotted vertical
lines drawn at each time when $\Tamb(t)$ crosses 0\degC\ for easy
identification of freeze and thaw events.

One sugar maple simulation (S1) is performed with ambient temperature
$\Tamb(t)$ taken equal to the smoothed temperature from
Figure~\ref{fig:tree-data}b-i and stem radius set to
$\Rtree=30.5$~cm. Because this is a more mature tree than the red
maples, we assume that the heartwood extends half-way through the stem
and take $\heartfrac=\frac{1}{2}$.  The resulting pressure curves are
displayed in Figure~\ref{fig:sims}c, with the experimental pressure
taken from the north-side sensor measurements (refer to
Figure~\ref{fig:tree-data}b-ii).

\section{Discussion}
\label{sec:discussion}

\subsection{Thaw events trigger pressure spikes}
\label{sec:discuss-data}

Based on the experimental measurements of pressure in the right-hand
plots of Figure~\ref{fig:tree-data}(right), a characteristic feature of
all three trees is the spikes or rapid increases in pressure that occur
at certain discrete times.  By viewing these plots alongside the
corresponding temperature curves from Figure~\ref{fig:tree-data}(left), 
it is clear that the spikes coincide with times that thaw events occur,
which is when the temperature increases past the freezing point (this
correspondence is much more clearly seen in the plots appearing later in
Figure~\ref{fig:sims}). After thawing, and for as long as the
temperature remains above zero, each spike is followed by a period
of gradually falling pressure where the rate of decrease appears fairly
consistent between thawing events.  This behavior is consistent with
other experimental pressure measurements in maple and related species
\citep{tyree-1983, cortes-sinclair-1985, ameglio-etal-2001,
  ewers-etal-2001}. The timing, amplitude and decay of these pressure
spikes will form the main points of comparison when we discuss the
numerical simulations in Section~\ref{sec:discuss-sims}.

It is also worth noting that for the sugar maple only, we have provided
two pressure curves (see Figure~\ref{fig:tree-data}b-ii) that correspond
to measurements taken from two sensors located opposite each other on
the north and south sides of the stem. The measured pressure variations
are qualitatively very similar, especially as regards the timing and
height of the pressure spikes and their subsequent decay.  A notable
exception is the two spikes recorded by the south sensor around day 20
that do not have matching spikes in the north side data. Because these
two spikes correspond to especially weak thaw events (``weak'' in the
sense that temperature exceeds zero for only a brief time interval) it
is likely that they are not experienced as thawing events around the
entire tree stem.  For this reason, we have chosen to use the north-side
pressure data for comparison with the S1 simulation.

\subsection{Xylem pressure is most sensitive to $\Rtree$ and $\gamma_s$} 
\label{sec:discuss-sensitivity}

The parameter sensitivity study in Figure~\ref{fig:sensitivity} will
guide our choice of parameter values for use in the experimental
comparisons in the next section, as well as pinpointing those parameters
that have the greatest impact on the model solution and hence are most
important to estimate accurately.  We begin by comparing the qualitative
features of the solutions in Figure~\ref{fig:sensitivity}a for different
stem radius, $\Rtree$. In all cases we observe that average stem
pressure $\overline{p}$ behaves similar to the experiments in that it
exhibits a steep increase whenever temperature increases past the
freezing point. This is followed by a gradual decline while the
temperature remains positive but as soon as temperature falls below zero
there is a similarly steep drop in pressure, after which the cycle
repeats. The amplitude of the pressure oscillations decreases for larger
radius trees, which is consistent with there being a larger sapwood area
to freeze and hence a correspondingly larger water uptake.  \mydel{What
  is perhaps more important for exudation is the ``envelope curves''
  that trace out max/min points of the pressure (see blue dashed curves
  in Figure~\ref{fig:sensitivity}a, left) which exhibit an upward trend
  that reflects a build-up in exudation pressure over time.}  The
pressure \mydel{envelope} is affected significantly by changes in
$\Rtree$ in that the amplitude of pressure oscillations decreases with
increased $\Rtree$; however, the time-averaged pressure (which sits
midway between the \mydel{envelope curves} \mymod{pressure
  minima/maxima}) is relatively insensitive to changes in radius.  This
contrasts with the relatively huge increase in root water uptake for
higher values of $\Rtree$, but again these larger trees have a
proportionally greater volume to freeze and over which to distribute
the stored pressure.

The simulations in Figure~\ref{fig:sensitivity}b correspond to different
values of the heartwood fraction $\heartfrac$ between 0 and 0.7.
Clearly, both the shape of the pressure oscillations and the exudation
pressure build-up are relatively insensitive to heartwood ratio. The
root water uptake curves shift downward as $\heartfrac$ increases, but
this is simply a geometric effect due to the corresponding decrease in
sapwood area. Among the two geometric parameters ($\Rtree$ and
$\heartfrac$) the solution is clearly most sensitive to stem radius,
which is fortunate since $\heartfrac$ is difficult to measure
(non-destructively) in a live tree.

The root area $\Atree$ that controls root water flux 
can be estimated from data in the literature, but it still remains the
most uncertain parameter in our model.  When increasing $\Atree$ over
two orders of magnitude (from 1 to 100), Figure~\ref{fig:sensitivity}c
shows that the amplitude of the stem pressure oscillations increases
significantly owing to a drop in the minimum pressure, whereas the
peak pressures remain approximately the same. However, over the
narrower range of values $\Atree\in[10.4,18.6]$ reported by
\citet{day-harris-2007} for red maples, the differences in stem
pressure are relatively small.

\leavethisout{The root reflection coefficient has a more conspicuous
  impact on exudation and Figure~\ref{fig:sensitivity}d shows that the
  oscillations in pressure trend upwards as $\Cout$ is reduced. This an
  obvious consequence of restricting root outflow which retains more
  water in the xylem that can freeze and subsequently store pressure in
  the fibers.  It is clearly important that we have incorporated the
  root reflection coefficient in our modified model equations and
  despite the relative insensitivity to this parameter, there is still
  an opportunity here for new experiments that aim to obtain better
  estimates of $\Cout$ in maple trees.}


The final parameter we consider is sap sugar concentration for which
Figure~\ref{fig:sensitivity}e shows that increasing $\gamma_s$ within
the range $[0.015, 0.04]$ has the greatest impact on increasing stem
pressure compared to the other three parameters.  This supports
our earlier remarks regarding the essential role played by sugar in
terms of generating a local differential between freeze/thaw in fibers
and vessels (due to FPD) which permits ice to accumulate in fibers at
the same time as the sap in neighboring vessels remains in the liquid
state. Increasing $\gamma_s$ therefore permits additional ice
accumulation in fibers through cryostatic suction which is also
reflected in a corresponding increase in root water uptake with
$\gamma_s$. It is worth recalling that sap sugar also induces an osmotic
contribution to pressure through Eq.~\eqref{eq:U-darcy}, but we have
clearly demonstrated in previous work~\citep{graf-ceseri-stockie-2015}
that osmosis is eclipsed in importance by the effect of FPD.  Finally,
because sap sugar content is so easy to measure, it is especially
important that any similar experimental study of sap pressure and
temperature also includes measurements of $\gamma_s$ in order that the
model can be properly calibrated.

\subsection{Multiscale model reproduces realistic pressure variations}
\label{sec:discuss-sims}

We begin by discussing the results in Figure~\ref{fig:sims}a,b that
compare experimental and numerical results for the two red maples R1 and
R2.  The measurements are dominated by pressure spikes appearing at
times $t\approx 13$, 15, 16, 18, 22 days that coincide with similar
spikes in simulations. Each spike is clearly matched with a thaw event
in which air temperature increases past \mymod{the critical 
  temperature}.  The simulated peak pressure for some spikes does not
reach the same peak value as in experiments but the correspondence is
nonetheless excellent, not to mention that the pressure minima following
subsequent freeze events are captured very closely. The simulations also
show that each spike is followed by a relaxation period during which the
pressure gradually decays at a rate that is similar to that seen in the
measured data.

During the initial 13~days on the other hand, the match between red
maple simulations and experiments is not nearly as close. The measured
R1 data in Figure~\ref{fig:sims}a shows that the pressure remains
essentially constant at zero whereas the simulation exhibits significant
pressure fluctuations in response to freeze--thaw events, most notably
on days 2 and 9. Deviations are also present with R2, although the match
is slighly better because the pressure sensor captures some small
fluctuations between days 1--2 and 8--11. One possible explanation for
these discrepancies is that the thaw events for times $t<13$~days are
weaker in that temperature increases only slightly above
0\degC\ before either falling below freezing again or hovering near
zero, which may be causing the stem to remain more deeply
frozen. These ``weak thaw'' events seem to be captured more readily by 
the model computations, although the computed pressure does exhibit a more
gradual increase instead of the sharp spike seen at the onset of other
thaw events; this behavior is especially apparent for the thaw event 
highlighted in Figure~\ref{fig:sims}a,b during days 9--11.  

The observations of low pressure/flow earlier in the season may also be
attributed to the fact that trees tend to be very well-buffered to
temperature on the north side due to reduced sun exposure
\citep{reid-driller-watson-2020} (recalling that our comparison is based
on north-side sensor data).  Furthermore, higher accumulations of snow
near the base of the stem may also limit the root water uptake.
In either case, it is well known that several repeated cycles of freeze
and thaw are typically required before pressure and flow rates can ramp up
to higher values. 

Next we shift to the comparison depicted in Figure~\ref{fig:sims}c for
the sugar maple S1, which shows that the timing of pressure spikes and
the subsequent relaxation rate from simulations both exhibit a
reasonable match with experiments, although the simulated pressure peak
values are significantly lower. We have as yet no definitive explanation
for this discrepancy, but it may be at least partly due to our
pressure average $\overline{p}$, which is integrated all the way to the
heartwood boundary and so includes portions of the xylem that may be
more deeply frozen and lower the average.  Finally, we single out the
thaw event on day 7 during which the simulated pressure builds up much
more gradually than in the experimental data, which is again analogous
to what we observed for the weak thaw events in trees R1 and R2.

\subsection{Essential mechanisms governing sap exudation}
\label{sec:conclusion}

These comparisons between experiments and simulations demonstrate that a
purely physical model is capable of capturing both qualitatively and
quantitatively the essential features of sap transport and pressure
generation observed in actual maple trees undergoing exudation.
Furthermore, we have clearly identified four mechanisms that are
essential for generating stem pressure build-up:
\begin{enumerate}
\item The \emph{distinctive cellular structure of maple sapwood} 
  which is made up of libriform fibers containing mostly gas that are
  connected hydraulically through selectively permeable walls to sap-filled
  vessels. This structure has two very important consequences: first, it
  provides a mechanism for fiber--vessel pressure exchange via
  compression of gas in the fibers; and second, the
  selectively-permeable nature of the fiber--vessel wall ensures that
  any liquid drawn into the fiber (through cryostatic suction) contains
  no sucrose.
  
\item The \emph{sugar content of sap} \mymod{affects the freeze--thaw
    induced pressure generation in maple xylem in two major
    ways}. First, dissolved sugar in the vessels induces a significant
  osmotic potential between fiber and vessel, which extends the range of
  pressures over which gas bubbles persist in sap.  \mymod{A second and
    more critical contribution} to exudation is the \emph{freezing point
    depression} (FPD) of roughly 0.16\degC\ in the vessel sap relative
  to pure water (assuming a 3\%\ sugar content). This is what allows the
  fibers to accumulate a frozen pure-ice layer while the sugary sap in
  neighboring vessels remains thawed because of the FPD. \mymod{The
    combination of these two effects is clearly evident in the
    comparison of pressure variations from the sensitivity study
    pictured in Figure~\ref{fig:sensitivity}a.}
  
\item \emph{A clear separation of spatial scales} that exists between
  freeze--thaw processes on the microscopic (cellular) scale and heat
  transport on the macroscopic (tree) scale.  Specifically, the FPD may
  appear to be insignificant on the macroscale on which freezing/thawing
  fronts propagate through the tree stem, but it dominates on the
  cellular scale by permitting thawed vessels to co-exist adjacent to
  partially frozen fibers.
\item An available supply of \emph{soil water in the liquid phase} (even
  under freezing conditions) is the key to generating a
  build-up in exudation pressure over multiple freeze--thaw cycles via
  accumulation of ice within the fibers. Although this feature is not
  specific to trees that exude, the availability of significant soil
  water under sub-zero conditions has recently been confirmed in
  experiments on maple saplings.
\end{enumerate}
Each of these distinguishing features has been recognized in other
studies of maple or related species; however, this is the first time
that they have all been linked together to construct a complete
quantitative model for the exudation process that also provides a reasonable
match with experimental measurements. This is a minimally complete model
in the sense that leaving out any of these four effects from the
governing equations results in a failure of the model tree stem to
accumulate exudation pressures that are consistent with actual trees.
It is important to recognize that these fundamental insights we
have gained into the physical mechanisms driving the sap exudation
process were only possible by developing a detailed mathematical model
and performing careful parametric studies of the resulting numerical
simulations.

\subsection{Opportunities for future research}
\label{sec:future}

The results presented in this paper and the mathematical model on
which they are based open up several opportunities for future research
in the study of exudation and sap flow in maple and other species.
\begin{itemize}
\item First of all, our sensitivity study singles out two parameters
  -- sap sugar content and root surface area -- for which accurate
  estimates are required in order to properly calibrate the
  model. Because the root surface area has not been accurately
  measured yet for maple trees, this suggests at least one opportunity
  for future experimental studies.
  
\item The 2D stem model can be extended in a straightforward fashion
  to a more realistic 3D axisymmetric stem geometry by stacking a
  sequence of 2D cross-sections in the vertical direction and then
  coupling sap flux and temperature between adjacent sections. We
  could then incorporate the effect of changes in gravitational
  potential with height, while at the same time obtaining a more
  realistic representation of how soil water from the roots is drawn
  by cryostatic suction to higher elevations in the tree.
    
\item With such a 3D model in hand, we could also easily incorporate
  the variation of sap sugar content with height that has been
  observed by \citet{milburn-zimmermann-1986}. It would then be
  natural to incorporate the temporal dynamics of the starch
  conversion process wherein living xylem cells release sugar into the
  vessels in response to temperature variations
  \citep{ameglio-etal-2001, wong-baggett-rye-2003}.
    
\item \mymod{Some authors have argued that hysteresis and supercooling
    play an important role in the freeze-thaw events that occur during
    the sap exudation process \citep{bozonnet-etal-2023, tyree-1983}.
    These effects would be relatively straightforward to incorporate
    into our model, but would have to be validated through careful
    comparisons with experiments on maple or other tree species
    \citep{charrier-etal-2015, neuner-xu-hacker-2010,
      robson-petty-1987}.}
  
\item Several studies suggest that the ability of trees to exude sap
  is can be attributed at least partially to root-level processes, in
  addition to the stem processes we have considered
  here~\citep{holtta-etal-2018, kramer-boyer-1995,
    westhoff-etal-2008}. A modelling study of root pressure generation
  would thus form a very interesting avenue for future research,
  especially in species such as birch which are thought to be
  dominated by root pressure.
  
\item Finally, there are fascinating connections to explore between sap
  exudation and freeze-induced embolism, motivated by studies that have
  demonstrated a close relationship between \mydel{between} embolism
  recovery and positive pressures in xylem~\citep{holtta-etal-2018,
    schenk-jansen-holtta-2021, sperry-etal-1988}.
\end{itemize}



%
%

\section*{Funding}

Natural Sciences and Engineering Research Council of Canada
(RGPIN-2021-04088 to JMS);  
Alexander von Humboldt Foundation (Fyodor Lynen Fellowship to IK); 
North American Maple Syrup Council Research Fund (to JMS); 
University of Vermont Agricultural Experiment Station (to TRW, 
TDP, AvdB).

\section*{\mymod{Data availability statement}}

\mymod{The experimental data and the Matlab source code used to perform
  the numerical simulations in this paper can be obtained from the
  corresponding author upon request.}

%

\section*{Authors' Contributions}

Study conception and design (JMS); mathematical model development (IK,
JMS); algorithm implementation and numerical simulations (IK, MZ, JMS);
experimental design and data collection (TRW, TDP, AvdB); data analysis,
synthesis and interpretation (MZ, TDP, AvdB, JMS); manuscript writing
and revision (all authors).





\begin{thebibliography}{56}
\expandafter\ifx\csname natexlab\endcsname\relax\def\natexlab#1{#1}\fi
\expandafter\ifx\csname url\endcsname\relax
  \def\url#1{\texttt{#1}}\fi
\expandafter\ifx\csname urlprefix\endcsname\relax\def\urlprefix{URL }\fi

\bibitem[{Allaire(1992)}]{allaire-1992}
Allaire G (1992). Homogenization and two-scale convergence. SIAM Journal on
  Mathematical Analysis 23(6):1482--1518.

\bibitem[{Am{\'e}glio et~al.(2001)Am{\'e}glio, Ewers, Cochard, Martignac,
  Vandame, Bodet, \& Cruiziat}]{ameglio-etal-2001}
Am{\'e}glio T, Ewers FW, Cochard H, Martignac M, Vandame M, Bodet C, Cruiziat P
  (2001). Winter stem xylem pressure in walnut trees: effects of carbohydrates,
  cooling and freezing. Tree Physiology 21(6):387--394.

\bibitem[{Baral et~al.(2017)Baral, Berninger, Schneider, \&
  Pothier}]{baral-etal-2017}
Baral SK, Berninger F, Schneider R, Pothier D (2017). Effects of heartwood
  formation on sugar maple (\emph{Acer saccharum} {M}arshall) discoloured wood
  proportion. Trees 31:105--114.

\bibitem[{Bozonnet et~al.(2023)Bozonnet, Saudreau, Badel, Am\'eglio, \&
  Charrier}]{bozonnet-etal-2023}
\mymod{Bozonnet C, Saudreau M, Badel E, Am\'eglio T, Charrier G (2023). Freeze
  dehydration vs.\ supercooling in tree stems: {P}hysical and physiological
  modelling. Tree Physiology (on-line),
  DOI:10.1093/treephys/tpad117.}

\bibitem[{Brodersen et~al.(2010)Brodersen, {McElrone}, Choat, Matthews, \&
    Shackel}]{brodersen-etal-2010}
  Brodersen CR, {McElrone} AJ, Choat B, Matthews MA, Shackel KA (2010). The
  dynamics of embolism repair in xylem: {I}n vivo visualizations using
  high-resolution computed tomography. Plant Physiology 154(3):1088--1095.

\bibitem[{Ceseri \& Stockie(2013)}]{ceseri-stockie-2013}
Ceseri M, Stockie JM (2013). A mathematical model for sap exudation in maple
  trees governed by ice melting, gas dissolution and osmosis. SIAM Journal on
  Applied Mathematics 73(2):649--676.

\bibitem[{Ceseri \& Stockie(2014)}]{ceseri-stockie-2014}
Ceseri M, Stockie JM (2014). A three-phase free boundary problem involving ice
  melting and gas dissolution. European Journal of Applied Mathematics
  25(4):449--480.

\bibitem[{Charrier et~al.(2015)Charrier, Pramsohler, Charra-Vaskou,
    Saudreau, Am{\'e}glio, Neuner, \& Mayr}]{charrier-etal-2015}
  \mymod{Charrier G, Pramsohler M, Charra-Vaskou K, Saudreau M,
    Am{\'e}glio T, Neuner G, Mayr S (2015). Ultrasonic emissions during
    ice nucleation and propagation in plant xylem. New Phytologist
    207:570--578.}
  
\bibitem[{Chavarr{\'i}a-Krauser \&
  Ptashnyk(2013)}]{chavarriakrauser-ptashnyk-2013}
Chavarr{\'i}a-Krauser A, Ptashnyk M (2013). Homogenization approach to water
  transport in plant tissues with periodic microstructures. Mathematical
  Modelling of Natural Phenomena 8(4):80--111.

\bibitem[{Cirelli et~al.(2008)Cirelli, Jagels, \& Tyree}]{cirelli-etal-2008}
Cirelli D, Jagels R, Tyree MT (2008). Toward an improved model of maple sap
  exudation: the location and role of osmotic barriers in sugar maple,
  butternut and white birch. Tree Physiology 28:1145--1155.

\bibitem[{Cortes \& Sinclair(1985)}]{cortes-sinclair-1985}
Cortes PM, Sinclair TR (1985). The role of osmotic potential in spring sap flow
  of mature sugar maple trees (\emph{Acer saccharum} {M}arsh.). Journal of
  Experimental Botany 36(1):12--24.


\bibitem[{Day \& Harris(2007)}]{day-harris-2007}
Day SD, Harris JR (2007). Fertilization of red maple (\emph{Acer rubrum}) and
  littleleaf linden (\emph{Tilia cordata}) trees at recommended rates does not
  aid tree establishment. Arboriculture and Urban Forestry 33(2):113--121.

\bibitem[{Duchesne et~al.(2016)Duchesne, Vincent, Wang, Ung, \&
  Swift}]{duchesne-etal-2016}
Duchesne I, Vincent M, Wang XA, Ung CH, Swift DE (2016). Wood mechanical
  properties and discoloured heartwood proportion in sugar maple and yellow
  birch grown in {N}ew {B}runswick. BioResources 11(1):2007--2019.

\bibitem[{Ewers et~al.(2001)Ewers, Am{\'e}glio, Cochard, Beaujard, Martignac,
  Vandame, Bodet, \& Cruiziat}]{ewers-etal-2001}
Ewers FW, Am{\'e}glio T, Cochard H, Beaujard F, Martignac M, Vandame M, Bodet
  C, Cruiziat P (2001). Seasonal variation in xylem pressure of walnut trees:
  root and stem pressures. Tree Physiology 21:1123--1132.

\bibitem[{Fowler \& Krantz(1994)}]{fowler-krantz-1994}
Fowler AC, Krantz WB (1994). A generalized secondary frost heave model. SIAM
  Journal on Applied Mathematics 54(6):1650--1675.

\bibitem[{Graf et~al.(2015)Graf, Ceseri, \& Stockie}]{graf-ceseri-stockie-2015}
Graf I, Ceseri M, Stockie JM (2015). Multiscale model of a freeze--thaw process
  for tree sap exudation. Journal of the Royal Society Interface 12:{20150665}.

\bibitem[{Graf \& Stockie(2014)}]{graf-stockie-2014}
Graf I, Stockie JM (Dec. 2014). A mathematical model for maple sap exudation.
  Maple Syrup Digest 26A(4):15--19.


\bibitem[{H{\"o}ltt{\"a} et~al.(2018)H{\"o}ltt{\"a}, {Dominguez Carrasco},
  Salmon, Aalto, Vanhatalo, B{\"a}ck, \& Lintunen}]{holtta-etal-2018}
H{\"o}ltt{\"a} T, {Dominguez Carrasco} MDR, Salmon Y, Aalto J, Vanhatalo A,
  B{\"a}ck J, Lintunen A (2018). Water relations in silver birch during
  springtime. {H}ow is sap pressurized? Plant Biology 20:834--847.


\bibitem[{Johnson(1945)}]{johnson-1945}
Johnson LPV (1945). Physiological studies on sap flow in the sugar maple,
  \emph{Acer saccharum} {M}arsh. Canadian Journal of Research C Botanical
  Sciences 23:192--197.

\bibitem[{Johnson \& Tyree(1992)}]{johnson-tyree-1992}
Johnson RW, Tyree MT (1992). Effect of stem water content on sap flow from
  dormant maple and butternut stems: {I}nduction of sap flow in butternut.
  Plant Physiology 100:853--858.

\bibitem[{Johnson et~al.(1987)Johnson, Tyree, \&
  Dixon}]{johnson-tyree-dixon-1987}
Johnson RW, Tyree MT, Dixon MA (1987). A requirement for sucrose in xylem sap
  flow from dormant maple trees. Plant Physiology 84:495--500.

\bibitem[{Jones \& Alli(1987)}]{jones-alli-1987}
Jones ARC, Alli I (1987). Sap yields, sugar content, and soluble carbohydrates
  of saps and syrups of some {C}anadian birch and maple species. Canadian
  Journal of Forest Research 17:263--266.


\bibitem[{{Konrad} et~al.(2017){Konrad}, Peter, \&
  Stockie}]{konrad-peter-stockie-2017}
{Konrad} I, Peter MA, Stockie JM (2017). A two-scale {S}tefan problem arising
  in a model for tree sap exudation. IMA Journal of Applied Mathematics
  82(4):726--762.

\bibitem[{Kramer \& Boyer(1995)}]{kramer-boyer-1995}
  Kramer PJ, Boyer JS (1995). The absorption of water and root and stem
  pressures. In: Water Relations of Plants and Soils. Academic Press, London,
  Ch.~6, pp. 167--200.
  
\bibitem[{Larochelle et~al.(1998)Larochelle, Forget, Rainville, \&
  Bousquet}]{larochelle-etal-1998}
Larochelle F, Forget {\'E}, Rainville A, Bousquet J (1998). Sources of temporal
  variation in sap sugar content in a mature sugar maple \emph{(Acer
  saccharum)} plantation. Forest Ecology and Management 106:307--313.

\bibitem[{Marvin(1949)}]{marvin-1949}
Marvin JW (1949). Changes in bark thickness during sap flow in sugar maples.
  Science 109(2827):231--232.

\bibitem[{Marvin(1958)}]{marvin-1958}
Marvin JW (1958). The physiology of maple sap flow. In: Thimann KV, Critchfield
  WB, Zimmermann MH (Eds.), The Physiology of Forest Trees: A Symposium Held at
  the Harvard Forest. Ronald Press, New York, pp. 95--124.

\bibitem[{Marvin(1968)}]{marvin-1968}
Marvin JW (Aug. 20-22, 1968). Physiology of sap production. In: Sugar Maple
  Conference. Michigan Technological University, Houghton, MI, pp. 12--15.

\bibitem[{Marvin \& Greene(1959)}]{marvin-greene-1959}
Marvin JW, Greene MT (Feb. 1959). Some factors affecting the yield from maple
  tapholes. Bulletin 611, Vermont Agricultural Experiment Station, Burlington,
  VT, available online at
  \url{http://cdi.uvm.edu/collections/item/bulletin611}.

\bibitem[{MATLAB(2020)}]{matlab-2020a}
MATLAB (2020). Version 9.8.0 (R2020a). The MathWorks Inc., Natick,
  Massachusetts.

\bibitem[{Merwin \& Lyon(1909)}]{merwin-lyon-1909}
Merwin HE, Lyon H (1909). Sap pressure in the birch stem. Botanical Gazette
  48(6):442--458.

\bibitem[{Milburn \& O'Malley(1984)}]{milburn-omalley-1984}
Milburn JA, O'Malley PER (1984). Freeze-induced sap absorption in \emph{Acer
  pseudoplatanus}: a possible mechanism. Canadian Journal of Botany
  62(10):2101--2106.

\bibitem[{Milburn \& Zimmermann(1986)}]{milburn-zimmermann-1986}
Milburn JA, Zimmermann MH (1986). Sapflow in the sugar maple in the leafless
  state. Journal of Plant Physiology 124:331--344.

\bibitem[{Neuner et~al.(2010)Neuner, Xu, \& Hacker}]{neuner-xu-hacker-2010}
\mymod{Neuner G, Xu B, Hacker J (2010). Velocity and pattern of ice propagation and
  deep supercooling in woody stems of \emph{Castanea sativa}, \emph{Morus
  nigra} and \emph{Quercus robur} measured by {IDTA}. Tree Physiology
  30:1037--1045.}

  
\bibitem[{Perkins \& {van den Berg}(2009)}]{perkins-vandenberg-2009}
Perkins TD, {van den Berg} AK (2009). Maple syrup---{P}roduction, composition,
  chemistry, and sensory characteristics. In: Taylor SL (Ed.), Advances in Food
  and Nutrition Research. Vol.~56. Elsevier, Ch.~4, pp. 101--143.

\bibitem[{Reid et~al.(2020)Reid, Driller, \& Watson}]{reid-driller-watson-2020}
Reid S, Driller T, Watson M (2020). A two-dimensional heat transfer model for
  predicting freeze-thaw events in sugar maple trees. Agricultural and Forest
  Meteorology 294:{108139}.

\bibitem[{Robitaille et~al.(1995)Robitaille, Boutin, \&
  Lachance}]{robitaille-boutin-lachance-1995}
Robitaille G, Boutin R, Lachance D (1995). Effects of soil freezing stress on
  sap flow and sugar content of mature sugar maples (\emph{Acer saccharum}).
  Canadian Journal of Forest Research 25(4):577--587.

\bibitem[{Robson \& Petty(1987)}]{robson-petty-1987}
\mymod{Robson DJ, Petty JA (1987). Freezing in conifer xylem: {I}. {P}ressure changes
  and growth velocity of ice. Journal of Experimental Botany 38(11):1901--1908.}

\bibitem[{Sachs(1860)}]{sachs-1860}
Sachs J (20 July 1860). Quellungserscheinungen an {H}\"olzern. Botanische
  Zeitung 18(29):253--259.

\bibitem[{Schenk et~al.(2021)Schenk, Jansen, \&
  H{\"o}ltt{\"a}}]{schenk-jansen-holtta-2021}
Schenk HJ, Jansen S, H{\"o}ltt{\"a} T (2021). Positive pressure in xylem and
  its role in hydraulic function. New Phytologist 230:27--45.

\bibitem[{Sorkin(2014)}]{sorkin-2014}
Sorkin L (Jan. 20, 2014). Maple syrup revolution: {A} new discovery could
  change the business forever. Modern Farmer,
  \url{http://modernfarmer.com/2014/01/maple-syrup-revolution}.

\bibitem[{Sperry et~al.(1988)Sperry, Donnelly, \& Tyree}]{sperry-etal-1988}
Sperry JS, Donnelly JR, Tyree MT (1988). Seasonal occurrence of xylem embolism
  in sugar maple (\emph{{A}cer saccharum}). American Journal of Botany
  75(8):1212--1218.


\bibitem[{Steudle \& Peterson(1998)}]{steudle-peterson-1998}
  \mymod{Steudle E, Peterson CA (1998). How does water get through roots? Journal of
    Experimental Biology 49(322):775--788.}
  
\bibitem[{Stevens \& Eggert(1945)}]{stevens-eggert-1945}
Stevens CL, Eggert RL (1945). Observations on the causes of the flow of sap in
  red maple. Plant Physiology 20:636--648.

\bibitem[{Tyree(1983)}]{tyree-1983}
Tyree MT (1983). Maple sap uptake, exudation, and pressure changes correlated
  with freezing exotherms and thawing endotherms. Plant Physiology 73:277--285.

\bibitem[{Tyree(1995)}]{tyree-1995}
Tyree MT (April 10-12, 1995). The mechanism of maple sap exudation. In:
  Terazawa M, McLeod CA, Tamai Y (Eds.), Tree Sap: Proceedings of the First
  International Symposium on Sap Utilization. Hokkaido University Press,
  Bifuka, Japan, pp. 37--45.


\bibitem[{Tyree \& Zimmermann(2002)}]{tyree-zimmermann-2002}
Tyree MT, Zimmermann MH (2002). Xylem Structure and the Ascent of Sap, 2nd
  Edition. Springer Series in Wood Science. Springer-Verlag, Berlin.

\bibitem[{Visintin(1996)}]{visintin-1996}
Visintin A (1996). Models of Phase Transitions. Vol.~28 of Progress in
  Nonlinear Differential Equations and Their Applications. Birkh\"auser,
  Boston.

\bibitem[{Westhoff et~al.(2008)Westhoff, Schneider, Zimmermann, Mimietz,
    Stinzing, Wegner, Kaiser, Krohne, Shirley, Jakob, Bamberg, Bentrup,
    \& Zimmermann}]{westhoff-etal-2008} Westhoff M, Schneider H,
    Zimmermann D, Mimietz S, Stinzing A, Wegner LH, Kaiser W, Krohne G,
    Shirley S, Jakob P, Bamberg E, Bentrup FW, Zimmermann U (2008).  The
    mechanisms of refilling of xylem conduits and bleeding of tall birch
    during spring. Plant Biology 10:604--623.

\bibitem[{Wiegand(1906)}]{wiegand-1906}
  Wiegand KM (1906). Pressure and flow of sap in wood. The American Naturalist
  40(474):409--453.

\bibitem[{Wilmot(2006)}]{wilmot-2006}
Wilmot T (Jun. 2006). Temperatures in the sugarbush. Maple Syrup Digest
  18A(2):20--23.

\bibitem[{Wong et~al.(2003)Wong, Baggett, \& Rye}]{wong-baggett-rye-2003}
Wong BL, Baggett KL, Rye AH (2003). Seasonal patterns of reserve and soluble
  carbohydrates in mature sugar maple (\emph{Acer saccharum}). Canadian Journal
  of Botany 81:780--788.


\bibitem[{Yang \& Tyree(1992)}]{yang-tyree-1992}
Yang S, Tyree MT (1992). A theoretical model of hydraulic conductivity recovery
  from embolism with comparison to experimental data on \emph{Acer saccharum}.
  Plant Cell Environment 15:633--643.

\bibitem[{Zarrinderakht(2017)}]{zarrinderakht-mscthesis-2017}
Zarrinderakht M (Mar. 2017). Numerical simulations of a multiscale model for
  maple sap exudation. Master's thesis, Department of Mathematics, Simon Fraser
  University, Burnaby, BC, Canada, available at
  \url{https://theses.lib.sfu.ca/thesis/etd10022}.

\bibitem[{Zwieniecki \& Holbrook(2000)}]{zwieniecki-holbrook-2000}
  Zwieniecki MA, Holbrook NM (2000). Bordered pit structure and vessel wall
  surface properties: {I}mplications for embolism repair. Plant Physiology
  123(3):1015--1020.
\end{thebibliography}



\providecommand{\noopsort}[1]{}

\end{document}